\documentclass{amsart}
\usepackage{anysize}

\newcommand{\tlowername}[2]%
{$\stackrel{\makebox[1pt]{#1}}%
{\begin{picture}(0,0)%
\put(0,0){\makebox(0,6)[t]{\makebox[1pt]{$#2$}}}%
\end{picture}}$}%

\newcommand{\AR}[1]%
{\begin{picture}(#1,0)%
\put(0,0){\vector(1,0){#1}}%
\end{picture}}%

\newcommand{\DOTAR}[1]%
{\NUMBEROFDOTS=#1%
\divide\NUMBEROFDOTS by 3%
\begin{picture}(#1,0)%
\multiput(0,0)(3,0){\NUMBEROFDOTS}{\circle*{1}}%
\put(#1,0){\vector(1,0){0}}%
\end{picture}}%

\newcommand{\MONO}[1]%
{\begin{picture}(#1,0)%
\put(0,0){\vector(1,0){#1}}%
\put(2,-2){\line(0,1){4}}%
\end{picture}}%

\newcommand{\EPI}[1]%
{\begin{picture}(#1,0)(-#1,0)%
\put(-#1,0){\vector(1,0){#1}}%
\put(-6,-2){\line(0,1){4}}%
\end{picture}}%

\newcommand{\BIMO}[1]%
{\begin{picture}(#1,0)(-#1,0)%
\put(-#1,0){\vector(1,0){#1}}%
\put(-6,-2){\line(0,1){4}}%
\put(-#1,-2){\hspace{2pt}\line(0,1){4}}%
\end{picture}}%

\newcommand{\BIAR}[1]%
{\begin{picture}(#1,4)%
\put(0,0){\vector(1,0){#1}}%
\put(0,4){\vector(1,0){#1}}%
\end{picture}}%

\newcommand{\EQL}[1]%
{\begin{picture}(#1,0)%
\put(0,1){\line(1,0){#1}}%
\put(0,-1){\line(1,0){#1}}%
\end{picture}}%

\newcommand{\ADJAR}[1]%
{\begin{picture}(#1,4)%
\put(0,0){\vector(1,0){#1}}%
\put(#1,4){\vector(-1,0){#1}}%
\end{picture}}%

\newcommand{\BKAR}[1]%
{\begin{picture}(#1,0)%
\put(#1,0){\vector(-1,0){#1}}%
\end{picture}}%

\newcommand{\BKDOTAR}[1]%
{\NUMBEROFDOTS=#1%
\divide\NUMBEROFDOTS by 3%
\begin{picture}(#1,0)%
\multiput(#1,0)(-3,0){\NUMBEROFDOTS}{\circle*{1}}%
\put(0,0){\vector(-1,0){0}}%
\end{picture}}%

\newcommand{\BKMONO}[1]%
{\begin{picture}(#1,0)(-#1,0)%
\put(0,0){\vector(-1,0){#1}}%
\put(-2,-2){\line(0,1){4}}%
\end{picture}}%

\newcommand{\BKEPI}[1]%
{\begin{picture}(#1,0)%
\put(#1,0){\vector(-1,0){#1}}%
\put(6,-2){\line(0,1){4}}%
\end{picture}}%

\newcommand{\BKBIMO}[1]%
{\begin{picture}(#1,0)%
\put(#1,0){\vector(-1,0){#1}}%
\put(6,-2){\line(0,1){4}}%
\put(#1,-2){\hspace{-2pt}\line(0,1){4}}%
\end{picture}}%

\newcommand{\BKBIAR}[1]%
{\begin{picture}(#1,4)%
\put(#1,0){\vector(-1,0){#1}}%
\put(#1,4){\vector(-1,0){#1}}%
\end{picture}}%

\newcommand{\BKADJAR}[1]%
{\begin{picture}(#1,4)%
\put(0,4){\vector(1,0){#1}}%
\put(#1,0){\vector(-1,0){#1}}%
\end{picture}}%

\newcommand{\lowername}[2]%
{$\stackrel{\makebox[1pt]{#1}}%
{\begin{picture}(0,0)%
\truex{600}%
\put(0,0){\makebox(0,\value{x})[t]{\makebox[1pt]{$#2$}}}%
\end{picture}}$}%

\newcommand{\hcase}[2]%
{\makebox[0pt]%
{\raisebox{-1pt}[0pt][0pt]{#1{#2}}}}%

\newcommand{\Hcase}[3]%
{\makebox[0pt]
{\raisebox{-1pt}[0pt][0pt]%
{$\stackrel{\makebox[0pt]{$\textstyle{#2}$}}{#1{#3}}$}}}%

\newcommand{\hcasE}[3]%
{\makebox[0pt]%
{\raisebox{-9pt}[0pt][0pt]%
{\lowername{#1{#3}}{#2}}}}%

\newcommand{\hbicase}[2]%
{\makebox[0pt]%
{\raisebox{-2.5pt}[0pt][0pt]{#1{#2}}}}%

\newcommand{\Hbicase}[4]%
{\makebox[0pt]
{\raisebox{-10.5pt}[0pt][0pt]%
{$\stackrel{\makebox[0pt]{$\textstyle{#2}$}}%
{\mbox{\lowername{#1{#4}}{#3}}}$}}}%

\newcommand{\EAR}[1]%
{\begin{picture}(#1,0)%
\put(0,0){\vector(1,0){#1}}%
\end{picture}}%

\newcommand{\EDOTAR}[1]%
{\truex{100}\truey{300}%
\NUMBEROFDOTS=#1%
\divide\NUMBEROFDOTS by \value{y}%
\begin{picture}(#1,0)%
\multiput(0,0)(\value{y},0){\NUMBEROFDOTS}%
{\circle*{\value{x}}}%
\put(#1,0){\vector(1,0){0}}%
\end{picture}}%

\newcommand{\EMONO}[1]%
{\begin{picture}(#1,0)%
\put(0,0){\vector(1,0){#1}}%
\truex{300}\truey{600}%
\put(\value{x},-\value{x}){\line(0,1){\value{y}}}%
\end{picture}}%

\newcommand{\EEPI}[1]%
{\begin{picture}(#1,0)(-#1,0)%
\put(-#1,0){\vector(1,0){#1}}%
\truex{300}\truey{600}\truez{800}%
\put(-\value{z},-\value{x}){\line(0,1){\value{y}}}%
\end{picture}}%

\newcommand{\EBIMO}[1]%
{\begin{picture}(#1,0)(-#1,0)%
\put(-#1,0){\vector(1,0){#1}}%
\truex{300}\truey{600}\truez{800}%
\put(-\value{z},-\value{x}){\line(0,1){\value{y}}}%
\put(-#1,-\value{x}){\hspace{3pt}\line(0,1){\value{y}}}%
\end{picture}}%

\newcommand{\EBIAR}[1]%
{\truex{400}%
\begin{picture}(#1,\value{x})%
\put(0,0){\vector(1,0){#1}}%
\put(0,\value{x}){\vector(1,0){#1}}%
\end{picture}}%

\newcommand{\EEQL}[1]%
{\begin{picture}(#1,0)%
\truex{200}%
\put(0,\value{x}){\line(1,0){#1}}%
\put(0,0){\line(1,0){#1}}%
\end{picture}}%

\newcommand{\EADJAR}[1]%
{\truex{400}%
\begin{picture}(#1,\value{x})%
\put(0,0){\vector(1,0){#1}}%
\put(#1,\value{x}){\vector(-1,0){#1}}%
\end{picture}}%

\newcommand{\ear}%
{\hspace{\SOURCE\unitlength}%
\hcase{\EAR}{\ARROWLENGTH}}%

\newcommand{\Ear}[1]%
{\hspace{\SOURCE\unitlength}%
\Hcase{\EAR}{#1}{\ARROWLENGTH}}%

\newcommand{\eaR}[1]%
{\hspace{\SOURCE\unitlength}%
\hcasE{\EAR}{#1}{\ARROWLENGTH}}%

\newcommand{\edotar}%
{\hspace{\SOURCE\unitlength}%
\hcase{\EDOTAR}{\ARROWLENGTH}}%

\newcommand{\Edotar}[1]%
{\hspace{\SOURCE\unitlength}%
\Hcase{\EDOTAR}{#1}{\ARROWLENGTH}}%

\newcommand{\edotaR}[1]%
{\hspace{\SOURCE\unitlength}%
\hcasE{\EDOTAR}{#1}{\ARROWLENGTH}}%

\newcommand{\emono}%
{\hspace{\SOURCE\unitlength}%
\hcase{\EMONO}{\ARROWLENGTH}}%

\newcommand{\Emono}[1]%
{\hspace{\SOURCE\unitlength}%
\Hcase{\EMONO}{#1}{\ARROWLENGTH}}%

\newcommand{\emonO}[1]%
{\hspace{\SOURCE\unitlength}%
\hcasE{\EMONO}{#1}{\ARROWLENGTH}}%

\newcommand{\eepi}%
{\hspace{\SOURCE\unitlength}%
\hcase{\EEPI}{\ARROWLENGTH}}%

\newcommand{\Eepi}[1]%
{\hspace{\SOURCE\unitlength}%
\Hcase{\EEPI}{#1}{\ARROWLENGTH}}%

\newcommand{\eepI}[1]%
{\hspace{\SOURCE\unitlength}%
\hcasE{\EEPI}{#1}{\ARROWLENGTH}}%

\newcommand{\ebimo}%
{\hspace{\SOURCE\unitlength}%
\hcase{\EBIMO}{\ARROWLENGTH}}%

\newcommand{\Ebimo}[1]%
{\hspace{\SOURCE\unitlength}%
\Hcase{\EBIMO}{#1}{\ARROWLENGTH}}%

\newcommand{\ebimO}[1]%
{\hspace{\SOURCE\unitlength}%
\hcasE{\EBIMO}{#1}{\ARROWLENGTH}}%

\newcommand{\eiso}%
{\hspace{\SOURCE\unitlength}%
\Hcase{\EAR}{\cong}{\ARROWLENGTH}}%

\newcommand{\Eiso}[1]%
{\hspace{\SOURCE\unitlength}%
\Hcase{\EAR}{\cong#1}{\ARROWLENGTH}}%

\newcommand{\eisO}[1]%
{\hspace{\SOURCE\unitlength}%
\hcasE{\EAR}{\cong#1}{\ARROWLENGTH}}%

\newcommand{\ebiar}%
{\hspace{\SOURCE\unitlength}%
\hbicase{\EBIAR}{\ARROWLENGTH}}%

\newcommand{\Ebiar}[2]%
{\hspace{\SOURCE\unitlength}%
\Hbicase{\EBIAR}{#1}{#2}{\ARROWLENGTH}}%

\newcommand{\eeql}%
{\hspace{\SOURCE\unitlength}%
\hbicase{\EEQL}{\ARROWLENGTH}}%

\newcommand{\eadjar}%
{\hspace{\SOURCE\unitlength}%
\hbicase{\EADJAR}{\ARROWLENGTH}}%

\newcommand{\Eadjar}[2]%
{\hspace{\SOURCE\unitlength}%
\Hbicase{\EADJAR}{#1}{#2}{\ARROWLENGTH}}%

\newcommand{\WAR}[1]%
{\begin{picture}(#1,0)%
\put(#1,0){\vector(-1,0){#1}}%
\end{picture}}%

\newcommand{\WDOTAR}[1]%
{\truex{100}\truey{300}%
\NUMBEROFDOTS=#1%
\divide\NUMBEROFDOTS by \value{y}%
\begin{picture}(#1,0)%
\multiput(#1,0)(-\value{y},0){\NUMBEROFDOTS}%
{\circle*{\value{x}}}%
\put(0,0){\vector(-1,0){0}}%
\end{picture}}%

\newcommand{\WMONO}[1]%
{\begin{picture}(#1,0)(-#1,0)%
\put(0,0){\vector(-1,0){#1}}%
\truex{300}\truey{600}%
\put(-\value{x},-\value{x}){\line(0,1){\value{y}}}%
\end{picture}}%

\newcommand{\WEPI}[1]%
{\begin{picture}(#1,0)%
\put(#1,0){\vector(-1,0){#1}}%
\truex{300}\truey{600}\truez{800}%
\put(\value{z},-\value{x}){\line(0,1){\value{y}}}%
\end{picture}}%

\newcommand{\WBIMO}[1]%
{\begin{picture}(#1,0)%
\put(#1,0){\vector(-1,0){#1}}%
\truex{300}\truey{600}\truez{800}%
\put(\value{z},-\value{x}){\line(0,1){\value{y}}}%
\put(#1,-\value{x}){\hspace{-3pt}\line(0,1){\value{y}}}%
\end{picture}}%

\newcommand{\WBIAR}[1]%
{\truex{400}%
\begin{picture}(#1,\value{x})%
\put(#1,0){\vector(-1,0){#1}}%
\put(#1,\value{x}){\vector(-1,0){#1}}%
\end{picture}}%

\newcommand{\WADJAR}[1]%
{\truex{400}%
\begin{picture}(#1,\value{x})%
\put(0,\value{x}){\vector(1,0){#1}}%
\put(#1,0){\vector(-1,0){#1}}%
\end{picture}}%

\newcommand{\war}%
{\hspace{\SOURCE\unitlength}%
\hcase{\WAR}{\ARROWLENGTH}}%

\newcommand{\War}[1]%
{\hspace{\SOURCE\unitlength}%
\Hcase{\WAR}{#1}{\ARROWLENGTH}}%

\newcommand{\waR}[1]%
{\hspace{\SOURCE\unitlength}%
\hcasE{\WAR}{#1}{\ARROWLENGTH}}%

\newcommand{\wdotar}%
{\hspace{\SOURCE\unitlength}%
\hcase{\WDOTAR}{\ARROWLENGTH}}%

\newcommand{\Wdotar}[1]%
{\hspace{\SOURCE\unitlength}%
\Hcase{\WDOTAR}{#1}{\ARROWLENGTH}}%

\newcommand{\wdotaR}[1]%
{\hspace{\SOURCE\unitlength}%
\hcasE{\WDOTAR}{#1}{\ARROWLENGTH}}%

\newcommand{\wmono}%
{\hspace{\SOURCE\unitlength}%
\hcase{\WMONO}{\ARROWLENGTH}}%

\newcommand{\Wmono}[1]%
{\hspace{\SOURCE\unitlength}%
\Hcase{\WMONO}{#1}{\ARROWLENGTH}}%

\newcommand{\wmonO}[1]%
{\hspace{\SOURCE\unitlength}%
\hcasE{\WMONO}{#1}{\ARROWLENGTH}}%

\newcommand{\wepi}%
{\hspace{\SOURCE\unitlength}%
\hcase{\WEPI}{\ARROWLENGTH}}%

\newcommand{\Wepi}[1]%
{\hspace{\SOURCE\unitlength}%
\Hcase{\WEPI}{#1}{\ARROWLENGTH}}%

\newcommand{\wepI}[1]%
{\hspace{\SOURCE\unitlength}%
\hcasE{\WEPI}{#1}{\ARROWLENGTH}}%

\newcommand{\wbimo}%
{\hspace{\SOURCE\unitlength}%
\hcase{\WBIMO}{\ARROWLENGTH}}%

\newcommand{\Wbimo}[1]%
{\hspace{\SOURCE\unitlength}%
\Hcase{\WBIMO}{#1}{\ARROWLENGTH}}%

\newcommand{\wbimO}[1]%
{\hspace{\SOURCE\unitlength}%
\hcasE{\WBIMO}{#1}{\ARROWLENGTH}}%

\newcommand{\wiso}%
{\hspace{\SOURCE\unitlength}%
\Hcase{\WAR}{\cong}{\ARROWLENGTH}}%

\newcommand{\Wiso}[1]%
{\hspace{\SOURCE\unitlength}%
\Hcase{\WAR}{#1}{\ARROWLENGTH}}%

\newcommand{\wisO}[1]%
{\hspace{\SOURCE\unitlength}%
\hcasE{\WAR}{#1}{\ARROWLENGTH}}%

\newcommand{\wbiar}%
{\hspace{\SOURCE\unitlength}%
\hbicase{\WBIAR}{\ARROWLENGTH}}%

\newcommand{\Wbiar}[2]%
{\hspace{\SOURCE\unitlength}%
\Hbicase{\WBIAR}{#1}{#2}{\ARROWLENGTH}}%

\newcommand{\weql}%
{\hspace{\SOURCE\unitlength}%
\hbicase{\EEQL}{\ARROWLENGTH}}%

\newcommand{\wadjar}%
{\hspace{\SOURCE\unitlength}%
\hbicase{\WADJAR}{\ARROWLENGTH}}%

\newcommand{\Wadjar}[2]%
{\hspace{\SOURCE\unitlength}%
\Hbicase{\WADJAR}{#1}{#2}{\ARROWLENGTH}}%

\newcommand{\vcase}[2]{#1{#2}}%

\newcommand{\Vcase}[3]{\makebox[0pt]%
{\makebox[0pt][r]{\raisebox{0pt}[0pt][0pt]{${#2}\hspace{2pt}$}}}#1{#3}}%

\newcommand{\vcasE}[3]{\makebox[0pt]%
{#1{#3}\makebox[0pt][l]{\raisebox{0pt}[0pt][0pt]{\hspace{2pt}$#2$}}}}%

\newcommand{\vbicase}[2]{\makebox[0pt]{{#1{#2}}}}%

\newcommand{\Vbicase}[4]{\makebox[0pt]%
{\makebox[0pt][r]{\raisebox{0pt}[0pt][0pt]{$#2$\hspace{4pt}}}#1{#4}%
\makebox[0pt][l]{\raisebox{0pt}[0pt][0pt]{\hspace{5pt}$#3$}}}}%

\newcommand{\SAR}[1]%
{\begin{picture}(0,0)%
\put(0,0){\makebox(0,0)%
{\begin{picture}(0,#1)%
\put(0,#1){\vector(0,-1){#1}}%
\end{picture}}}\end{picture}}%

\newcommand{\SDOTAR}[1]%
{\truex{100}\truey{300}%
\NUMBEROFDOTS=#1%
\divide\NUMBEROFDOTS by \value{y}%
\begin{picture}(0,0)%
\put(0,0){\makebox(0,0)%
{\begin{picture}(0,#1)%
\multiput(0,#1)(0,-\value{y}){\NUMBEROFDOTS}%
{\circle*{\value{x}}}%
\put(0,0){\vector(0,-1){0}}%
\end{picture}}}\end{picture}}%

\newcommand{\SMONO}[1]%
{\begin{picture}(0,0)%
\put(0,0){\makebox(0,0)%
{\begin{picture}(0,#1)%
\put(0,#1){\vector(0,-1){#1}}%
\truex{300}\truey{600}%
\put(0,#1){\begin{picture}(0,0)%
\put(-\value{x},-\value{x}){\line(1,0){\value{y}}}\end{picture}}%
\end{picture}}}\end{picture}}%

\newcommand{\SEPI}[1]%
{\begin{picture}(0,0)%
\put(0,0){\makebox(0,0)%
{\begin{picture}(0,#1)%
\put(0,#1){\vector(0,-1){#1}}%
\truex{300}\truey{600}\truez{800}%
\put(-\value{x},\value{z}){\line(1,0){\value{y}}}%
\end{picture}}}\end{picture}}%

\newcommand{\SBIMO}[1]%
{\begin{picture}(0,0)%
\put(0,0){\makebox(0,0)%
{\begin{picture}(0,#1)%
\put(0,#1){\vector(0,-1){#1}}%
\truex{300}\truey{600}\truez{800}%
\put(0,#1){\begin{picture}(0,0)%
\put(-\value{x},-\value{x}){\line(1,0){\value{y}}}\end{picture}}%
\put(-\value{x},\value{z}){\line(1,0){\value{y}}}%
\end{picture}}}\end{picture}}%

\newcommand{\SBIAR}[1]%
{\begin{picture}(0,0)%
\truex{200}%
\put(0,0){\makebox(0,0)%
{\begin{picture}(0,#1)\put(-\value{x},#1){\vector(0,-1){#1}}%
\put(\value{x},#1){\vector(0,-1){#1}}%
\end{picture}}}\end{picture}}%

\newcommand{\SEQL}[1]%
{\begin{picture}(0,0)%
\truex{100}%
\put(0,0){\makebox(0,0)%
{\begin{picture}(0,#1)\put(-\value{x},#1){\line(0,-1){#1}}%
\put(\value{x},#1){\line(0,-1){#1}}%
\end{picture}}}\end{picture}}%

\newcommand{\sarv}[1]{\vcase{\SAR}{#100}}%

\newcommand{\sar}{\sarv{50}}%

\newcommand{\Sarv}[2]{\Vcase{\SAR}{#1}{#200}}%

\newcommand{\Sar}[1]{\Sarv{#1}{50}}%

\newcommand{\saRv}[2]{\vcasE{\SAR}{#1}{#200}}%

\newcommand{\saR}[1]{\saRv{#1}{50}}%

\newcommand{\Sisov}[2]%
{\Vbicase{\SAR}{#1\hspace{-2pt}}{\hspace{-2pt}\cong}{#200}}%

\newcommand{\seqlv}[1]{\vbicase{\SEQL}{#100}}%

\newcommand{\seql}{\seqlv{50}}%

\newcommand{\NAR}[1]%
{\begin{picture}(0,0)%
\put(0,0){\makebox(0,0)%
{\begin{picture}(0,#1)\put(0,0){\vector(0,1){#1}}%
\end{picture}}}\end{picture}}%

\newcommand{\NDOTAR}[1]%
{\truex{100}\truey{300}%
\NUMBEROFDOTS=#1%
\divide\NUMBEROFDOTS by \value{y}%
\begin{picture}(0,0)%
\put(0,0){\makebox(0,0)%
{\begin{picture}(0,#1)%
\multiput(0,0)(0,\value{y}){\NUMBEROFDOTS}%
{\circle*{\value{x}}}%
\put(0,#1){\vector(0,1){0}}%
\end{picture}}}\end{picture}}%

\newcommand{\NMONO}[1]%
{\begin{picture}(0,0)%
\put(0,0){\makebox(0,0)%
{\begin{picture}(0,#1)%
\put(0,0){\vector(0,1){#1}}%
\truex{300}\truey{600}%
\put(-\value{x},\value{x}){\line(1,0){\value{y}}}%
\end{picture}}}%
\end{picture}}%

\newcommand{\NEPI}[1]%
{\begin{picture}(0,0)%
\put(0,0){\makebox(0,0)%
{\begin{picture}(0,#1)%
\put(0,0){\vector(0,1){#1}}%
\truex{300}\truey{600}\truez{800}%
\put(0,#1){\begin{picture}(0,0)%
\put(-\value{x},-\value{z}){\line(1,0){\value{y}}}\end{picture}}%
\end{picture}}}\end{picture}}%

\newcommand{\NBIMO}[1]%
{\begin{picture}(0,0)%
\put(0,0){\makebox(0,0)%
{\begin{picture}(0,#1)%
\put(0,0){\vector(0,1){#1}}%
\truex{300}\truey{600}\truez{800}%
\put(-\value{x},\value{x}){\line(1,0){\value{y}}}%
\put(0,#1){\begin{picture}(0,0)%
\put(-\value{x},-\value{z}){\line(1,0){\value{y}}}\end{picture}}%
\end{picture}}}\end{picture}}%

\newcommand{\NBIAR}[1]%
{\begin{picture}(0,0)%
\truex{200}%
\put(0,0){\makebox(0,0)%
{\begin{picture}(0,#1)\put(-\value{x},0){\vector(0,1){#1}}%
\put(\value{x},0){\vector(0,1){#1}}%
\end{picture}}}\end{picture}}%

\newcommand{\Nisov}[2]%
{\Vbicase{\NAR}{#1\hspace{-2pt}}{\hspace{-2pt}\cong}{#200}}%

\newcommand{\NEDOTAR}%
{\truex{100}\truey{212}%
\NUMBEROFDOTS=5800%
\divide\NUMBEROFDOTS by \value{y}%
\begin{picture}(0,0)%
\multiput(-2900,-2900)(\value{y},\value{y}){\NUMBEROFDOTS}%
{\circle*{\value{x}}}%
\put(2900,2900){\vector(1,1){0}}%
\end{picture}}%

\newcommand{\SWDOTAR}%
{\truex{100}\truey{212}%
\NUMBEROFDOTS=5800%
\divide\NUMBEROFDOTS by \value{y}%
\begin{picture}(0,0)%
\multiput(2900,2900)(-\value{y},-\value{y}){\NUMBEROFDOTS}%
{\circle*{\value{x}}}%
\put(-2900,-2900){\vector(-1,-1){0}}%
\end{picture}}%

\newcommand{\SEDOTAR}%
{\truex{100}\truey{212}%
\NUMBEROFDOTS=5800%
\divide\NUMBEROFDOTS by \value{y}%
\begin{picture}(0,0)%
\multiput(-2900,2900)(\value{y},-\value{y}){\NUMBEROFDOTS}%
{\circle*{\value{x}}}%
\put(2900,-2900){\vector(1,-1){0}}%
\end{picture}}%

\newcommand{\NWDOTAR}%
{\truex{100}\truey{212}%
\NUMBEROFDOTS=5800%
\divide\NUMBEROFDOTS by \value{y}%
\begin{picture}(0,0)%
\multiput(2900,-2900)(-\value{y},\value{y}){\NUMBEROFDOTS}%
{\circle*{\value{x}}}%
\put(-2900,2900){\vector(-1,1){0}}%
\end{picture}}%

\newcommand{\ENEAR}[2]%
{\makebox[0pt]{\begin{picture}(0,0)%
\put(0,-150){\makebox(0,0){\begin{picture}(0,0)%
\put(-6600,-3300){\vector(2,1){13200}}%
\truex{200}\truey{800}\truez{600}%
\put(-\value{x},\value{x}){\makebox(0,\value{z})[r]{${#1}$}}%
\put(\value{x},-\value{y}){\makebox(0,\value{z})[l]{${#2}$}}%
\end{picture}}}\end{picture}}}%

\newcommand{\ESEAR}[2]%
{\makebox[0pt]{\begin{picture}(0,0)%
\put(0,-150){\makebox(0,0){\begin{picture}(0,0)%
\put(-6600,3300){\vector(2,-1){13200}}%
\truex{200}\truey{800}\truez{600}%
\put(\value{x},\value{x}){\makebox(0,\value{z})[l]{${#1}$}}%
\put(-\value{x},-\value{y}){\makebox(0,\value{z})[r]{${#2}$}}%
\end{picture}}}\end{picture}}}%

\newcommand{\WNWAR}[2]%
{\makebox[0pt]{\begin{picture}(0,0)%
\put(0,-150){\makebox(0,0){\begin{picture}(0,0)%
\put(6600,-3300){\vector(-2,1){13200}}%
\truex{200}\truey{800}\truez{600}%
\put(\value{x},\value{x}){\makebox(0,\value{z})[l]{${#1}$}}%
\put(-\value{x},-\value{y}){\makebox(0,\value{z})[r]{${#2}$}}%
\end{picture}}}\end{picture}}}%

\newcommand{\WSWAR}[2]%
{\makebox[0pt]{\begin{picture}(0,0)%
\put(0,-150){\makebox(0,0){\begin{picture}(0,0)%
\put(6600,3300){\vector(-2,-1){13200}}%
\truex{200}\truey{800}\truez{600}%
\put(-\value{x},\value{x}){\makebox(0,\value{z})[r]{${#1}$}}%
\put(\value{x},-\value{y}){\makebox(0,\value{z})[l]{${#2}$}}%
\end{picture}}}\end{picture}}}%

\newcommand{\NNEAR}[2]%
{\raisebox{-1pt}[0pt][0pt]{\begin{picture}(0,0)%
\put(0,0){\makebox(0,0){\begin{picture}(0,0)%
\put(-3300,-6600){\vector(1,2){6600}}%
\truex{100}\truez{600}%
\put(-\value{x},\value{x}){\makebox(0,\value{z})[r]{${#1}$}}%
\put(\value{x},-\value{z}){\makebox(0,\value{z})[l]{${#2}$}}%
\end{picture}}}\end{picture}}}%

\newcommand{\SSWAR}[2]%
{\raisebox{-1pt}[0pt][0pt]{\begin{picture}(0,0)%
\put(0,0){\makebox(0,0){\begin{picture}(0,0)%
\put(3300,6600){\vector(-1,-2){6600}}%
\truex{100}\truez{600}%
\put(-\value{x},\value{x}){\makebox(0,\value{z})[r]{${#1}$}}%
\put(\value{x},-\value{z}){\makebox(0,\value{z})[l]{${#2}$}}%
\end{picture}}}\end{picture}}}%

\newcommand{\SSEAR}[2]%
{\raisebox{-1pt}[0pt][0pt]{\begin{picture}(0,0)%
\put(0,0){\makebox(0,0){\begin{picture}(0,0)%
\put(-3300,6600){\vector(1,-2){6600}}%
\truex{200}\truez{600}%
\put(\value{x},\value{x}){\makebox(0,\value{z})[l]{${#1}$}}%
\put(-\value{x},-\value{z}){\makebox(0,\value{z})[r]{${#2}$}}%
\end{picture}}}\end{picture}}}%

\newcommand{\NNWAR}[2]%
{\raisebox{-1pt}[0pt][0pt]{\begin{picture}(0,0)%
\put(0,0){\makebox(0,0){\begin{picture}(0,0)%
\put(3300,-6600){\vector(-1,2){6600}}%
\truex{200}\truez{600}%
\put(\value{x},\value{x}){\makebox(0,\value{z})[l]{${#1}$}}%
\put(-\value{x},-\value{z}){\makebox(0,\value{z})[r]{${#2}$}}%
\end{picture}}}\end{picture}}}%

\newcommand{\Necurve}[2]%
{\begin{picture}(0,0)%
\truex{1300}\truey{2000}\truez{200}%
\put(0,\value{x}){\oval(#200,\value{y})[t]}%
\put(0,\value{x}){\makebox(0,0){\begin{picture}(#200,0)%
\put(#200,0){\vector(0,-1){\value{z}}}%
\put(0,0){\line(0,-1){\value{z}}}\end{picture}}}%
\truex{2500}%
\put(0,\value{x}){\makebox(0,0)[b]{${#1}$}}%
\end{picture}}%

\newcommand{\Nwcurve}[2]%
{\begin{picture}(0,0)%
\truex{1300}\truey{2000}\truez{200}%
\put(0,\value{x}){\oval(#200,\value{y})[t]}%
\put(0,\value{x}){\makebox(0,0){\begin{picture}(#200,0)%
\put(#200,0){\line(0,-1){\value{z}}}%
\put(0,0){\vector(0,-1){\value{z}}}\end{picture}}}%
\truex{2500}%
\put(0,\value{x}){\makebox(0,0)[b]{${#1}$}}%
\end{picture}}%

\newcommand{\Securve}[2]%
{\begin{picture}(0,0)%
\truex{1300}\truey{2000}\truez{200}%
\put(0,-\value{x}){\oval(#200,\value{y})[b]}%
\put(0,-\value{x}){\makebox(0,0){\begin{picture}(#200,0)%
\put(#200,0){\vector(0,1){\value{z}}}%
\put(0,0){\line(0,1){\value{z}}}\end{picture}}}%
\truex{2500}%
\put(0,-\value{x}){\makebox(0,0)[t]{${#1}$}}%
\end{picture}}%

\newcommand{\Swcurve}[2]%
{\begin{picture}(0,0)%
\truex{1300}\truey{2000}\truez{200}%
\put(0,-\value{x}){\oval(#200,\value{y})[b]}%
\put(0,-\value{x}){\makebox(0,0){\begin{picture}(#200,0)%
\put(#200,0){\line(0,1){\value{z}}}%
\put(0,0){\vector(0,1){\value{z}}}\end{picture}}}%
\truex{2500}%
\put(0,-\value{x}){\makebox(0,0)[t]{${#1}$}}%
\end{picture}}%

\newcommand{\Escurve}[2]%
{\begin{picture}(0,0)%
\truex{1400}\truey{2000}\truez{200}%
\put(\value{x},0){\oval(\value{y},#200)[r]}%
\put(\value{x},0){\makebox(0,0){\begin{picture}(0,#200)%
\put(0,0){\vector(-1,0){\value{z}}}%
\put(0,#200){\line(-1,0){\value{z}}}\end{picture}}}%
\truex{2500}%
\put(\value{x},0){\makebox(0,0)[l]{${#1}$}}%
\end{picture}}%

\newcommand{\Encurve}[2]%
{\begin{picture}(0,0)%
\truex{1400}\truey{2000}\truez{200}%
\put(\value{x},0){\oval(\value{y},#200)[r]}%
\put(\value{x},0){\makebox(0,0){\begin{picture}(0,#200)%
\put(0,0){\line(-1,0){\value{z}}}%
\put(0,#200){\vector(-1,0){\value{z}}}\end{picture}}}%
\truex{2500}%
\put(\value{x},0){\makebox(0,0)[l]{${#1}$}}%
\end{picture}}%

\newcommand{\Wscurve}[2]%
{\begin{picture}(0,0)%
\truex{1300}\truey{2000}\truez{200}%
\put(-\value{x},0){\oval(\value{y},#200)[l]}%
\put(-\value{x},0){\makebox(0,0){\begin{picture}(0,#200)%
\put(0,0){\vector(1,0){\value{z}}}%
\put(0,#200){\line(1,0){\value{z}}}\end{picture}}}%
\truex{2400}%
\put(-\value{x},0){\makebox(0,0)[r]{${#1}$}}%
\end{picture}}%

\newcommand{\Wncurve}[2]%
{\begin{picture}(0,0)%
\truex{1300}\truey{2000}\truez{200}%
\put(-\value{x},0){\oval(\value{y},#200)[l]}%
\put(-\value{x},0){\makebox(0,0){\begin{picture}(0,#200)%
\put(0,0){\line(1,0){\value{z}}}%
\put(0,#200){\vector(1,0){\value{z}}}\end{picture}}}%
\truex{2400}%
\put(-\value{x},0){\makebox(0,0)[r]{${#1}$}}%
\end{picture}}%

\newcount\SCALE%

\newcount\NUMBER%

\newcount\LINE%

\newcount\COLUMN%

\newcount\WIDTH%

\newcount\SOURCE%

\newcount\ARROW%

\newcount\TARGET%

\newcount\ARROWLENGTH%

\newcount\NUMBEROFDOTS%

\newcounter{x}%

\newcounter{y}%

\newcounter{z}%

\newcounter{horizontal}%

\newcounter{vertical}%

\newskip\itemlength%

\newskip\firstitem%

\newskip\seconditem%

\newcommand{\printarrow}{}%

\newcommand{\truex}[1]{%
\NUMBER=#1%
\multiply\NUMBER by 100%
\divide\NUMBER by \SCALE%
\setcounter{x}{\NUMBER}}%

\newcommand{\truey}[1]{%
\NUMBER=#1%
\multiply\NUMBER by 100%
\divide\NUMBER by \SCALE%
\setcounter{y}{\NUMBER}}%

\newcommand{\truez}[1]{%
\NUMBER=#1%
\multiply\NUMBER by 100%
\divide\NUMBER by \SCALE%
\setcounter{z}{\NUMBER}}%

\newcommand{\changecounters}[1]{%
\SOURCE=\ARROW%
\ARROW=\TARGET%
\settowidth{\itemlength}{#1}%
\ifdim \itemlength > 2800\unitlength%
\addtolength{\itemlength}{-2800\unitlength}%
\TARGET=\itemlength%
\divide\TARGET by 1310%
\multiply\TARGET by 100%
\divide\TARGET by \SCALE%
\else%
\TARGET=0%
\fi%
\ARROWLENGTH=5000%
\advance\ARROWLENGTH by -\SOURCE%
\advance\ARROWLENGTH by -\TARGET%
\advance\SOURCE by -\TARGET}%

\newcommand{\initialize}[1]{%
\LINE=0%
\COLUMN=0%
\WIDTH=0%
\ARROW=0%
\TARGET=0%
\changecounters{#1}%
\renewcommand{\printarrow}{#1}%
\begin{center}%
\vspace{10pt}%
\begin{picture}(0,0)}%

\newcommand{\DIAGV}[2]{%
\SCALE=#1%
\setlength{\unitlength}{655sp}%
\multiply\unitlength by \SCALE%
\divide\unitlength by 100%
\initialize{\mbox{$#2$}}}%

\newcommand{\n}[1]{%
\changecounters{\mbox{$#1$}}%
\put(\COLUMN,\LINE){\makebox(0,0){\printarrow}}%
\thinlines%
\renewcommand{\printarrow}{\mbox{$#1$}}%
\advance\COLUMN by 4000}%

\newcommand{\nn}[1]{%
\put(\COLUMN,\LINE){\makebox(0,0){\printarrow}}%
\thinlines%
\ifnum \WIDTH < \COLUMN%
\WIDTH=\COLUMN%
\else%
\fi%
\advance\LINE by -4000%
\COLUMN=0%
\ARROW=0%
\TARGET=0%
\changecounters{\mbox{$#1$}}%
\renewcommand{\printarrow}{\mbox{$#1$}}}%

\newcommand{\conclude}{%
\put(\COLUMN,\LINE){\makebox(0,0){\printarrow}}%
\thinlines%
\ifnum \WIDTH < \COLUMN%
\WIDTH=\COLUMN%
\else%
\fi%
\setcounter{horizontal}{\WIDTH}%
\setcounter{vertical}{-\LINE}%
\end{picture}}%

\newcommand{\diag}{%
\conclude%
\raisebox{0pt}[0pt][\value{vertical}\unitlength]{}%
\hspace*{\value{horizontal}\unitlength}%
\vspace{10pt}%
\end{center}%
\setlength{\unitlength}{1pt}}%

\newcommand{\diagv}[3]{%
\conclude%
\NUMBER=#1%
\rule{0pt}{\NUMBER pt}%
\hspace*{-#2pt}%
\raisebox{0pt}[0pt][\value{vertical}\unitlength]{}%
\hspace*{\value{horizontal}\unitlength}
\NUMBER=#3%
\advance\NUMBER by 10%
\vspace*{\NUMBER pt}%
\end{center}%
\setlength{\unitlength}{1pt}}%

\newcommand{\N}[1]%
{\raisebox{0pt}[7pt][0pt]{$#1$}}%

\newcommand{\crosslength}[2]{%
\settowidth{\firstitem}{#1}%
\settowidth{\seconditem}{#2}%
\ifdim\firstitem < \seconditem%
\itemlength=\seconditem%
\else%
\itemlength=\firstitem%
\fi%
\divide\itemlength by 2%
\hspace{\itemlength}}%

\marginsize{2cm}{2cm}{3cm}{3cm}
\newtheorem{thm}{Theorem}[section]
\newtheorem{cor}[thm]{Corollary}
\newtheorem{lem}[thm]{Lemma}
\newtheorem{prop}[thm]{Proposition}
\theoremstyle{definition}
\newtheorem{defn}[thm]{Definition}
\theoremstyle{remark}
\newtheorem{rem}[thm]{Remark}
\numberwithin{equation}{section}

\def\Mod{\mbox{-Mod}}

\newcommand{\Natur}{{\mathbb N}}

\newcommand{\Hom}{\mbox{\rm Hom}}

\newcommand{\Ext}{\mbox{\rm Ext}}

\def\O{{\mathcal O}}

\def\P{{\mathcal P}}
\def\dgP{{\it dg}\widetilde{{\mathcal P}}}
\def\barP{{\widetilde{{\mathcal P}}}}
\def\dgUo{{\it dg}\widetilde{{\mathcal U}^{\perp}}}
\def\barUo{{\widetilde{{\mathcal U}^{\perp}}}}
\def\Qco{\mathfrak{Qco}}
\def\fiz{\leftarrow}

\def\Z{\mathbb{Z}}
\newcommand{\U}{\mathcal{U}}
\def\Uo{{\mathcal U}^{\perp}}
\def\oUo{{^{\perp}}({\mathcal U}^{\perp})}
\newcommand{\Ch}{{\rm Ch}}
\def\li{{\displaystyle \lim_{\rightarrow}}\ }

\begin{document}
\baselineskip=15pt
\title[projective model structure]{Locally projective monoidal model structure for complexes of quasi-coherent sheaves on ${\bf P^1(k)}$}%

\author{E. Enochs, S. Estrada and J.R. Garc\'{\i}a-Rozas}
\address{Department of Mathematics, University of Kentucky, Lexington, Kentucky
40506-0027, U.S.A.} \email{enochs@ms.uky.edu}

\address{
Departamento de Matem\'atica Aplicada, Universidad de Murcia,Campus
del Espinardo, Espinardo (Murcia) 30100, Spain}
\email{sestrada@um.es}%
\address{
Departamento de \'Algebra y A. Matem\'atico, Universidad de
Almer\'{\i}a, Almer\'{\i}a 04071, Spain}\email{jrgrozas@ual.es}



\begin{abstract}
We will generalize the projective model structure in the category of
unbounded complexes of modules over a commutative ring to the
category of unbounded complexes of quasi-coherent sheaves over the
projective line. Concretely we will define a locally projective
model structure in the category of complexes of quasi-coherent
sheaves on the projective line. In this model structure the
cofibrant objects are the dg-locally projective complexes. We also
describe the fibrations of this model structure and show that the
model structure is monoidal. We point out that this model structure
is necessarily different from other known model structures such as
the injective model structure and the locally free model structure.
\end{abstract}
\maketitle

\section{Introduction}
Quasi-coherent sheaves are known to play the role of modules in
algebraic geometry. And even from a homological viewpoint they
behave much like modules. For example the derived category of
$R\Mod$, $\mathcal{D}(R)$ (here $R$ is a commutative ring with
identity) is well understood because there are several Quillen model
structures on $\Ch(R)$ (the category of unbounded complexes of
$R$-modules) which allow one to define and compute the extension
functors. These are the projective model structure and the injective
model structure. It is known that the injective model structure is
not suitable for studying the torsion functors since this structure,
is not compatible with the graded tensor product of $\Ch(R)$,
induced from the tensor product of $R\Mod$. But the projective model
structure is compatible with the tensor product (see \cite[Chapter
4]{Hov}) so it can be used to define the torsion functors.
Furthermore it has been recently proved, by using the (positive)
solution to the flat cover conjecture (cf. \cite{BiBash}), that
there is an induced flat model structure which is compatible with
the tensor product (see \cite{Gill}).

Now let us consider the category $\Qco(X)$ of quasi-coherent sheaves
over a scheme $X$. It has been proved in \cite{Est} that this is a
Grothendieck category, so hence we can apply a result due to Joyal
(in \cite{Beke}) that inherits an injective model structure which
allows one to compute the derived extension functors in the category
of quasi-coherent sheaves on any scheme. However there is a natural
tensor product in $\Qco(X)$, so it would be desirable to impose a
model structure in $\Ch(\Qco(X))$ compatible with the tensor product
of quasi-coherent sheaves. The main problem is that $\Qco(X)$ does
not have enough projectives, so the problem is of a different nature
from that of the case of $R$-modules. In some circumstances the
existence of a family of flat generators can be used to replace the
projective ones. For example in \cite{Gill2} is proved that the
category of unbounded complexes of sheaves of $\O$-modules admits a
flat model structure similar to that of $\Ch(R)$ by using the fact
that there are enough flat objects in the category of sheaves of
modules on a commutative ring. But it is not known in general if the
category of quasi-coherent sheaves on an arbitrary scheme admits a
family of flat generators. In \cite{Est} a family is exhibited that
makes $\Qco(X)$ into a locally $\lambda$-presentable category, for
$\lambda$ a certain regular cardinal. But they are not flat in
general. However for some nice schemes (which are in practise the
most used for algebraic geometers) like quasi-compact and
quasi-separated there are enough flat objects, so a modified version
of the results of \cite{Gill2} together with the positive solution
of the flat cover conjecture given in \cite{Est} allow one to impose
a flat model structure in $\Qco(X)$, at least in this case ($X$
quasi-compact and quasi-separated).

Now let us fix our scheme $X$ to be a closed immersion of the
projective space ${\bf P^n}(k)$ ($k$ is a field). Then there is a
nice family of generators for $\Qco(X)$ with finite projective
dimension. We have the family of $\O(m)$, $m\in\mathbb{Z}$ for ${\bf
P^n}(k)$. These give the family $\{i^*(\O(m)):\ m\in \mathbb{Z}\}$,
where $i:X\hookrightarrow {\bf P^n}(k)$ (see \cite[pg.
120]{Hartshorne} for notation and terminology) we will let $\O(m)$
denote $i^*(\O(m))$. Of course they are not projective but in some
circumstances they behave like projective objects. For instance, for
the case $n=1$ a classic result of Grothendieck states that every
finitely generated and free quasi-coherent sheaf decomposes as the
direct sum of $\O(m)$'s. Our goal in this paper will be to show that
this generators allow to get what we call a locally projective model
structure in $\Ch(\Qco({\bf P^1}(k))$ which is compatible with the
closed symmetrical monoidal structure of $\Qco({\bf P^1}(k))$. This
may be surprising at first since the class of locally projective
quasi-coherent sheaves is contained strictly in the class of flat
quasi-coherent sheaves and one could have the impression that for
categories without enough projectives but with enough flat objects,
the flat model structure would be the ``smallest'' one which is
compatible with the tensor product of the category.

The main idea we use to get our result is a generalized version of
Kaplansky's theorem (see \cite[Theorem 1]{Kap}) which states that
every locally projective quasi-coherent sheaf on ${\bf P^1}(k)$ is a
direct transfinite extension of countably generated quasi-coherent
sheaves (Theorem \ref{p1}). Direct (and inverse) transfinite
extensions are widely studied in \cite{EnIa}.

The paper is structured as follows: in Section \ref{locfree} we
introduce the cotorsion pair cogenerated by the class of locally
free generators when $X$ is a scheme with enough locally frees. In
Section \ref{locproj} we particularize the previous situation to the
scheme ${\bf P^1}(k)$ and we are able to prove that locally
projective quasi-coherent sheaves appear in the left side of a
cotorsion pair (Subsection \ref{coinc}). We also give a complete
description of the right side of this cotorsion pair in Subsection
\ref{desc}. Section \ref{tool} is devoted to developing the tools we
need in proving that we have an induced model structure in
$\Ch(\Qco({\bf P^1}(k)))$. And finally in Section \ref{fin} we get
the locally projective monoidal structure in $\Ch(\Qco({\bf
P^1}(k)))$. We note that since the complete description given in
Section \ref{locproj} of the quasi-coherent sheaves involved in the
cotorsion pair cogenerated by $\{\O(m):\ m\in \Z\}$ we are able to
describe the fibrations and the cofibrations in the locally
projective monoidal model structure. Thus we note that our model
structure is necessarily different from that in \cite[Theorem
2.4]{Hov2}.


\section{Preliminaries}

In this section we introduce all basic definitions we need along the
paper. We recall from \cite{EnIa} the definition of a direct
transfinite extension. Let ${\mathcal A}$ be a Grothendieck
category. A direct system $ ( A_{ \alpha } | \alpha \leq \lambda ) $
is said to be continuous if $ A_0  = 0 $ and if for each limit
ordinal $ \beta \leq \lambda $ we have $ \displaystyle
A_{\beta}=\lim_{ \rightarrow } \:  A_{ \alpha } \:$ with the limit
over the $ \alpha < \beta $. The direct system $ ( A_{ \alpha}
|\alpha \leq \lambda ) $ is said to be a system of monomorphisms if
all the morphisms in the system are monomorphisms.

\begin{defn}
 An object $A$ of $ {\mathcal A} $ is said to be a direct
transfinite extension of objects of $ {\mathcal L} $ (here
${\mathcal L}$ is a class of objects of ${\mathcal A}$ closed under
isomorphisms) if $ A = \displaystyle \lim_{ \rightarrow } \: A_{
\alpha } $ for a continuous direct system $ ( A_{ \alpha } |  \alpha
\leq \lambda ) $ of monomorphisms  such that $ \mbox{coker} \: ( A_{
\alpha } \rightarrow A_{ \alpha + 1 } )$ is in $ {\mathcal L} $
whenever $ \alpha + 1 \leq \lambda $. The class $ {\mathcal L} $ is
said to be closed under direct transfinite extensions if each direct
transfinite extension of objects in $ {\mathcal L} $ is also in $
{\mathcal L} $.
\end{defn}

\begin{defn}
Given a class $\mathcal D$ of objects of ${\mathcal{A}}$ then the
class of objects $Y$ of ${\mathcal{A}}$ such that $\Ext^1(D,Y)=0$
for all $D\in {\mathcal D}$ is denoted ${\mathcal D}^{\perp}$.
Similarly $^{\perp}{\mathcal D}$ denotes the class of $Z$ such that
$\Ext^1(Z,D)=0$ for all $D\in {\mathcal D}$.
\end{defn}

\begin{defn}
A pair (${\mathcal F}, {\mathcal C}$) of classes of objects of
${{\mathcal{A}}}$ is said to be cotorsion pair if ${\mathcal
F}^{\perp}={\mathcal C}$ and if $^{\perp}{ \mathcal C}={\mathcal
F}$. It is said to have enough injectives (resp. enough projectives)
if for each $Z$ of ${{\mathcal{A}}}$ there exists an exact sequence
 $0\rightarrow Y\rightarrow C\rightarrow
F\rightarrow 0$ (resp. an exact sequence $0\rightarrow C'\rightarrow
F'\rightarrow Z\rightarrow 0$) where $F,F'\in {\mathcal F}$ and
where $C,C'\in {\mathcal C}$. The cotorsion pair $({\mathcal
F},{\mathcal C})$ is said to be complete if it has enough injectives
and projectives for each object of $\mathcal{A}$. We furthermore say
that $({\mathcal F}, {\mathcal C}) $ is functorially complete if
these sequences can be chosen in a functorial manner (depending on
$Y$ and on $Z$) (see Definition 2.3 of \cite{hovey2}).
\end{defn}

We also recall from \cite{enochja} the definition of a Kaplansky
class.
\begin{defn}
A class $\mathcal L$ of objects of ${\mathcal{A}}$ is said to be a
Kaplansky class  of ${\mathcal{A}}$ if for each cardinal $\kappa$
there is a cardinal $\lambda$ such that if $S\subset L$ for some
$L\in \mathcal L$ where $|S|\leq \kappa$ then there is an $L'\subset
L$ with $S\subset L'$ where $|L'|\leq \lambda $ and where $L'$ and
$L/L'$ are both in $\mathcal L$.
\end{defn}

Now let $X$ be a scheme. In \cite{Est} is described a category
equivalent to the category $\Qco(X)$ of quasi-coherent sheaves on
$X$. If we let ${\mathcal T}$ to denote the collection of all affine
open subsets, then we can define a quiver $Q=(V,E)$, where the
vertices are the affine opens and there is an arrow $v\to w$
whenever $w\subset v$. Then we let $R$ be a functor from $Q$ to the
category of commutative rings, given by $R(v)={\mathcal O}_X(v)$,
where ${\mathcal O}_X$ is the structure sheaf of $X$. Then the
category of quasi-coherent $R$-modules (see \cite[Section 2]{Est})
over $Q$ is equivalent to the category of quasi-coherent sheaves on
$X$.

Now if $X$ is the projective line $X={\bf P^1}(k)$ on a field $k$,
we will give a explicit description of the quiver point of view of
$\Qco(X)$. The structure sheaf ${\mathcal O}$ of the scheme ${\bf
P^1}(k)=(Proj\ k[x_0,x_1], {\mathcal O})$ may be identified with the
representation
$${\mathcal O}\equiv k[x]\hookrightarrow k[x,x^{-1}]\hookleftarrow k[x^{-1}]$$ of the
quiver
$$Q\equiv\bullet\to\bullet\fiz\bullet,$$ by just making the change
$x=x_1/x_0$. Then, a quasi-coherent sheaf of modules over ${\bf
P^1}(k)$ is a representation
$$M\stackrel{f}{\to} P\stackrel{g}{\fiz} N$$ of $Q$ with $M\in
k[x]\Mod$, $N\in k[x^{-1}]$, $P\in k[x,x^{-1}]\Mod$, $f$ an
$k[x]$-linear map and $g$ an $k[x^{-1}]$-linear map, satisfying that
$$S^{-1} f:S^{-1} M\to S^{-1} P\cong P$$ and $$T^{-1} g:T^{-1} N\to T^{-1}
P\cong P$$ are $k[x,x^{-1}]$-isomorphisms, where
$S=\{1,x,x^2,\cdots\}$ and $T=\{1,x^{-1},x^{-2},\cdots\}$. It is
easy to see that kernels and cokernels in these categories of
representations of $Q$ with those relations can be computed
componentwise.

\section{The locally free cotorsion pair in
$\Qco(S)$}\label{locfree}

\begin{thm}\label{und}
Let $S$ be a noetherian scheme with enough locally frees. Let us
denote by $\U$ the set of all locally free generators. Then the pair
$(\oUo,\Uo)$ is a complete cotorsion pair. Furthermore every $X\in
\oUo$ is a locally projective quasi-coherent sheaf on $S$.
\end{thm}
\begin{proof} It is clear that $(\oUo,\Uo)$ is a cotorsion pair. Let
us see that it is a complete cotorsion pair. By \cite[Lemma
1]{EkTrl} it follows that $\oUo$ contains all direct transfinite
extensions of the locally frees $S\in \U$. Furthermore by
\cite{EkTrl} (the arguments there are for modules but easily carry
over to our setting) for all $M\in \Qco(X)$ there exists a short
exact sequence
$$0\to M\to Y\to Z\to 0$$ where $Y\in \Uo$ and $Z$ a direct
transfinite extension of $S\in \U$ (so $Z\in \oUo$ by the previous).
This shows that the cotorsion pair $(\oUo,\Uo)$ has enough
injectives. To show that it has enough projectives we mimic the
Salce trick (see \cite[Lemmas 2.2 and 2.3]{Salce}). Given any $M\in
\Qco(X)$, since $X$ has enough locally frees there exists a short
exact sequence
$$0\to U\to V\to M\to 0 $$ where $V$ is a direct sum of $S\in \U$.
Now let
$$0\to U\to Y\to Z\to 0$$ be exact with $Y\in \Uo$ and $Z$ a direct
transfinite extension of $S\in \U$. Form a pushout and get
\DIAGV{50} {} \n{} \n{0}\n{}\n{0}\nn
{}\n{}\n{\sar}\n{}\n{\sar}\nn{0}
\n{\ear}\n{U}\n{\ear}\n{V}\n{\ear}\n{M}\n{\ear}\n{0}\nn
{}\n{}\n{\sar}\n{}\n{\sar}\n{}\n{\saR{id_{M}}}\nn
{0}\n{\ear}\n{Y}\n{\ear}\n{W}\n{\ear}\n{M}\n{\ear}\n{0}\nn
{}\n{}\n{\sar}\n{}\n{\sar}\nn {}\n{}\n{Z}\n{\Ear{id_{Z}}}\n{Z}\nn
{}\n{}\n{\sar}\n{}\n{\sar}\nn {}\n{}\n{0}\n{}\n{0}\diag

Then since $V$ is a direct sum of $S\in \U$ and since $Z$ is a
direct transfinite extension of $S\in \U$ we see that $W$ is a
direct transfinite extension of locally frees. Also $Y\in \Uo$.
Hence if $M\in \oUo$ we get that $0\to Y\to W\to M\to 0$ splits and
so $M$ is a direct summand of a direct transfinite extension of
$S$'s. But then it follows that, for every vertex $v$, $M(v)$ is a
direct summand of a transfinite extension of projective modules
$S(v)$'s, so by \cite[Proposition 3]{aus} $M$ is a locally
projective quasi-coherent sheaf. \end{proof}

\begin{rem} In the next sections we exhibit a case where the converse of
the previous result is also true, that is $\oUo$ consists precisely
of locally projective quasi-coherent sheaves.
\end{rem}
\section{Locally projective quasi-coherent sheaves on ${\bf
P^1}(k)$}\label{locproj}

Henceforth we let our scheme to be the projective line on a field
$k$. In this case we will prove that the model structure induced by
the complete cotorsion pair $(\oUo,\Uo)$ is the generalization of
the usual projective model structure on $\Ch(R)$ (for $R$ any
commutative ring) (see \cite[Section 2.3]{Hov}). We hope the methods
of the next sections will apply to more general schemes to give more
information about we can call ``the locally projective model
structure'' on $\Ch(\Qco({\bf P^1}(k)))$.

We will start by describing the elements of the class $\Uo$.
\subsection{Computation of $\Ext^1(\O(n),M\to P\fiz N)$}\label{desc}

\bigskip\par\noindent

Let us take $M\stackrel{\sigma}{\to}P\stackrel{\tau}{\fiz}N$. We use
the Baer description of $\Ext^1$. Given a short exact sequence
\begin{equation}\label{eq1}
0\to (M\to P\fiz N)\to (A\to C\fiz B)\to (k[x]\hookrightarrow
k[x,x^{-1}]\hookleftarrow k[x^{-1}])\to 0
\end{equation}
 we know that $0\to M\to
A\to k[x]\to 0$ is exact. So it is split exact. So we can let
$A=M\oplus k[x]$ and assume $$0\to M\to M\oplus k[x]\to k[x]\to 0$$
is the obvious exact sequence. Likewise we can take $C=P\oplus
k[x,x^{-1}]$ and $B=N\oplus k[x^{-1}]$.

So $$(A\to C\fiz B)= M\oplus k[x]\to P\oplus k[x,x^{-1}]\fiz N\oplus
k[x^{-1}].$$ By the exact sequence $(\ref{eq1})$ we see that the map
$$M\oplus k[x]\to P\oplus k[x,x^{-1}]$$ is completely determined by
a map $k[x]\to P$, that is, that given $y\in P$, if we consider the
map $$(m,p(x))\mapsto (\sigma(m)+p(x)\ y,p(x))$$ from $M\oplus k[x]$
to $P\oplus k[x,x^{-1}]$ then if we localize at
$S=\{1,x,x^2,\cdots\}$ then we get that $$S^{-1}((M\oplus k[x])\to
(P\oplus k[x,x^{-1}]))$$ an isomorphism. This is because
$S^{-1}(M\to P)$ and $S^{-1}(k[x]\hookrightarrow k[x,x^{-1}]$ are
both isomorphisms, the exactness of $S^{-1}(.)$ and the snake lemma.

So using any $y\in P$ we get a commutative \DIAGV{90} {M}
\n{\Ear{\sigma}} \n{P} \nn {\sar}  \n{} \n{\sar} \nn {M\oplus k[x]}
\n{\ear} \n{P\oplus k[x,x^{-1}]} \diag

Similarly given a $z\in P$ we get a commutative \DIAGV{90} {P}
\n{\War{\tau}} \n{N} \nn {\sar}  \n{} \n{\sar} \nn {P\oplus
k[x,x^{-1}]} \n{\war} \n{N\oplus k[x^{-1}]} \diag where the bottom
map is $(n,q(x^{-1}))\mapsto (\tau(n)+q(x^{-1})z,q(x^{-1}))$.

So all the above gives.
\begin{prop}
Any extension of $\Ext^1(\O(0),M\to P\fiz N)$ is completely
determined by arbitrary $y,z\in P$.

\end{prop}

Using the same sort of reasoning we can see that a section for
$$0\to(M\to P\fiz N)\to (M\oplus k[x]\to P\oplus k[x,x^{-1}]\fiz
N\oplus k[x^{-1}])$$ $$\to (k[x]\hookrightarrow
k[x,x^{-1}]\hookleftarrow k[x^{-1}])\to 0$$ (where the central term
is determined by a $u\in M$, $v\in N$) so $k[x]\to M\oplus k[x]$
maps 1 to $(u,1)$ and $k[x^{-1}]\to N\oplus k[x^{-1}]$ maps 1 to
$(v,1)$. Here the conditions on $u,v$ in order that we have a
morphism
$$(k[x]\hookrightarrow k[x,x^{-1}]\hookleftarrow k[x^{-1}])\to
(M\oplus k[x]\to P\oplus k[x,x^{-1}]\fiz N\oplus k[x^{-1}])$$ are
that $\sigma(u)+y=\tau(v)+z$, or that $\sigma(u)-\tau(v)=z-y$. Note
that since $k[x]\cap k[x^{-1}]=k$ this condition is all that is
needed in order to have a morphism. Since $y,z\in P$ are arbitrary,
$z-y$ can be any element of $P$. So we have proved the following.
\begin{prop}
$\Ext^1(\O(0),M\stackrel{\sigma}{\to} P\stackrel{\tau}{\fiz} N)=0$
if, and only if, $P=\sigma(M)+\tau(N)$. In fact
$$\Ext^1(\O(0),M\stackrel{\sigma}{\to} P\stackrel{\tau}{\fiz} N)\cong
P/(\sigma(M)+\tau(N))$$
\end{prop}
Using the same type argument we can get.
\begin{prop}
For any integer $n$, $$\Ext^1(\O(n),M\stackrel{\sigma}{\to}
P\stackrel{\tau}{\fiz} N)\cong P/(x^n\sigma(M)+\tau(N))$$ so
$\Ext^1(\O(n),M\stackrel{\sigma}{\to} P\stackrel{\tau}{\fiz} N)=0 $
if, and only if, $x^n\sigma(M)+\tau(N)=P$.
\end{prop}
\subsection{The class $\oUo$ coincides with the class of locally
projectives}\label{coinc}

\bigskip\par\noindent

Let us denote by $\P$ the class of all locally projective
quasi-coherents sheaves on ${\bf P^1}(k)$. So $(M\to P\fiz N)\in \P$
if and only if $M,P$ and $N$ are projective $k[x],k[x,x^{-1}]$ and
$k[x^{-1}]$-modules respectively. By Theorem \ref{und} we already
know that $\oUo\subseteq \P$. We prove now that the converse is also
true. We note that this Theorem is a generalization of a Kaplansky's
theorem (\cite[Theorem 1]{Kap}) for quasi-coherent sheaves on ${\bf
P^1}(k)$.
\begin{thm}\label{p1}
Any locally projective $(M\to P\fiz N)\in \P$ is a direct
transfinite extension of countably generated locally projective
quasi-coherent sheaves on ${\bf P^1}(k)$.
\end{thm}

\begin{proof} Since $M\to P\fiz N$ is a quasi-coherent sheaf on ${\bf
P^1}(k)$ then $P\cong S^{-1}M\cong T^{-1}N$ with
$S=\{1,x,x^2,\cdots\}$, $T=\{1,x^{-1},x^{-2},\cdots\}$. Suppose that
$M$ and $N$ are projective. By \cite{Kap} $M=\oplus_{i\in I}M_i$ and
$N=\oplus_{j\in J}N_j$ with each $M_i$ and $N_j$ countably
generated.

Let $I'\subset I$ be any countable subset. Then it is clear that
$$S^{-1}(\oplus_{i\in I'}M_i)\subset T^{-1}(\oplus_{j\in J'}N_j)
$$for some countable subset $J'\subset J$. Then let $I'\subset
I''\subset I$, $I''$ countable be such that $$S^{-1}(\oplus_{i\in
I''}M_i)\supset T^{-1}(\oplus_{j\in J'}N_j).
$$ Then continuing this zig-zag procedure we construct $$I'\subset
I''\subset I'''\subset\cdots \subset I$$ $$J'\subset J''\subset
J'''\subset\cdots \subset J$$ with each of $I',I'',\cdots,
J',J'',\cdots$ countable and satisfying the obvious conditions. Then
if $\overline{I}=\cup_{n\geq 1}I^{(n)}$, $\overline{J}=\cup_{n\geq
1}J^{(n)}$ we get $S^{-1}(\oplus_{i\in I}M_i)=T(\oplus_{j\in J}N_j)$
i.e. we have the subrepresentation $$\oplus_{i\in
\overline{I}}M_i\to
S^{-1}(\oplus_{i\in\overline{I}})M_i=T^{-1}(\oplus_{j\in\overline{J}})N_j$$
with $\oplus_{i\in \overline{I}}M_i$ and
$(\oplus_{j\in\overline{J}})N_j $ countably generated projective
modules. Notice that the quotient of the original $M\to P\fiz N$ by
this subrepresentation is isomorphic to a representation
$$\oplus_{i\in I-\overline{I}}M_i\to U\fiz \oplus_{i\in
J-\overline{J}}N_j.$$ We repeat the procedure with this
representation and see that we can find $\overline{I}\subset
\overline{\overline{I}}\subset I$, $\overline{J}\subset
\overline{\overline{J}}\subset J$, $\overline{\overline{I}}$,
$\overline{\overline{J}}$ countable, with $S^{-1}(\oplus_{i\in
\overline{\overline{I}}}M_i)=T^{-1}(\oplus_{j\in
\overline{\overline{J}}}M_j)$ but where $\overline{\overline{J}}$
contains any given countable subset of $J-\overline{J}$. So we
continue this procedure and see that we can write
$I=\cup_{\alpha<\lambda}I_{\alpha}$,
$J=\cup_{\alpha<\lambda}J_{\alpha}$ as continuous unions of subsets
($\lambda$ some ordinal number) such that $S^{-1}(\oplus_{i\in
I_{\alpha}}M_i)=T^{-1}(\oplus_{j\in J_{\alpha}}N_j)$ for each
$\alpha$ and such that if $\alpha+1<\lambda$ then
$I_{\alpha+1}-I_{\alpha}$ and $J_{\alpha+1}-J_{\alpha}$ are
countable. \end{proof}

\medskip\par\noindent
\begin{rem}
If a module is a direct transfinite extension of countably generated
projective modules, then it is a direct sum of countably generated
projective modules and conversely. But in the sheaf situation above
we do not get such a direct sum.
\end{rem}
\begin{thm}\label{p2}
Every countably generated and locally free quasi-coherent sheaf on
${\bf P^1}(k)$ is a direct transfinite extension of $\O(n)$'s.
\end{thm}
\begin{proof} Let $M\to P\fiz N$ be such that $M$ and $N$ are free with
given countable bases. Then, as usual, we can assume
$M\stackrel{id}{\to}P\fiz N$ (see for example \cite{EnTor}). So
$M\to P\fiz N$ will be given by an infinite matrix $$\left(
                  \begin{array}{ccc}
                    p_{11} & p_{12} & \cdots \\
                    p_{21} &\ldots  &  \\
                    \vdots &  & \ddots \\
                  \end{array}
                \right)$$
where the columns correspond to the image of the base elements of
$N$. Hence we have a column finite matrix. As usual we can assume
that the matrix is in upper triangular form i.e. is equal to
$$\left(
                  \begin{array}{ccccc}
                    p_{11} & p_{12} & p_{13} & \cdots \\
                    0 & p_{22} & p_{23} & \cdots  &  \\
                     0 &0 & p_{33} &\cdots \\
                     \vdots&\vdots  &\vdots & \ddots \\
                  \end{array}
                \right)$$
This matrix corresponds to an automorphism of $k[x,x^{-1}]\oplus
k[x,x^{-1}]\oplus \cdots$ so has an inverse (also column finite). So
we get that in fact each $p_{ii}$ is a unit of $k[x,x^{-1}]$ (we
need that the inverse matrix is also upper triangular). So we can
assume $p_{ii}=x^{n_i}$ for each $i=1,2,\cdots$. Then we see that
$M\to P\fiz N$ has $\O(n_1)$ as a subrepresentation (generated by
the first base elements of $M$ and $N$) and that the quotient of
$M\to P\fiz N$ by this $\O(n_1)$ has $\O(n_2)$ as a
subrepresentation, etc. So we see that $M\to P\fiz N$ is in fact a
direct transfinite extension with the corresponding quotients equal
to the $\O(n_1), \O(n_2),\cdots$'s. Hence $M\to P\fiz N$ is a direct
transfinite extension of $\O(n)$'s.
\end{proof}
Combining Theorems \ref{p1} and \ref{p2} we get
\begin{cor}
Any locally projective sheaf is the direct transfinite extension of
$\O(n)$'s.
\end{cor}
\begin{rem}
It seems unlikely that we can get any kind of uniqueness result or
even that any such sheaf is a direct sum of $\O(n)$'s. So this
result supports the claim that it is worthwhile studying transfinite
extensions.
\end{rem}
\section{Complete cotorsion pairs in $\Ch(\Qco({\bf
P^1}(k)))$}\label{tool}

We devote the next two sections to proving that the class of locally
projective quasi-coherent sheaves on ${\bf P^1}(k)$ induces a
Quillen's model structure in the category $\Ch(\Qco({\bf P^1}(k)))$
of unbounded complexes of quasi-coherent sheaves on ${\bf P^1(k)}$.

This will be done by applying Hovey's criteria (see \cite[Theorem
2.2]{hovey2}) relating cotorsion pairs with model categories
structures. We need to recall some standard definitions concerning
to complexes. Most of them can found in \cite{jr}. We let
$(M,\delta)$ or simply $M$ to denote a complex
$$\cdots\to
M^{-1}\stackrel{\delta^{-1}}{\longrightarrow}
M^0\stackrel{\delta^{0}}{\longrightarrow}M^1\stackrel{\delta^{1}}{\longrightarrow}\cdots$$
If $M$ is a complex we let $Z(M)=\cdots \to Z_nM\to Z_{n+1}M\to
\cdots$ and $B(M)=\cdots \to B_nM\to B_{n+1}M\to \cdots$ be the
subcomplexes of cycles and boundaries of $M$.

For a given quasi-coherent sheaf $M$ we let $\underline{M}[n]$ the
complex with $M$ in the $-n$th place and $0$ in other places ($n\in
\Z$). We denote by $\overline{M}[n]$ the complex $\cdots\to 0\to
M\stackrel{id}{\to}M\to 0\to \cdots$ where the $M$ in the $-n-1$ and
$-n$th position ($n\in \Z)$.

If $(M,\delta_M)$ and $(N,\delta_N)$ are two chain complexes we
define $Hom(M,N)$ to be the complex
$$\cdots\to \prod_{k\in \Z}\Hom(M^k,N^{k+n})\stackrel{\delta^n}
{\to}\prod_{k\in \Z}\Hom(M^k,N^{k+n+1})\to \cdots,$$ where
$(\delta^n f)^k={\delta}_N^{k+n}f^k-(-1)^nf^{k+1}{\delta}_M^k$. Then
we define $\Ext_{\Ch(\Qco({\bf P^1}(k)))}(M,N)$ to be the group of
equivalence classes of short exact sequences of complexes $0\to N\to
L\to M\to 0$. We note that $\Ch(\Qco({\bf P^1}(k)))$ is a
Grothendieck category having the set
$\mathcal{I}=\{\underline{\O(m)}[n]:\ m,n\in \Z\}$ as a family of
generators. So $\Ext^i$, $i\in \Z$ functors can be computed by using
injective resolutions.

We recall from \cite{Gill} the following definitions: an exact
complex $E$ of quasi-coherent sheaves is said to be an $\Uo$ complex
if $E$ is exact and $Z_nE\in \Uo$, $\forall n\in \Z$. We let
$\barUo$ to denote the class of all $\Uo$ complexes. Then a complex
$M$ of quasi-coherent sheaves is said to be dg-locally projective if
$Hom(M,E)$ is an exact complex for any complex $E\in \barUo$. We let
$\dgP$ to denote the class of all dg-locally projective
quasi-coherent sheaves.

There are the corresponding dual definitions of the previous
classes. So we get the classes $\barP$ and $\dgUo$ of locally
projective complexes and dg-$\Uo$ complexes of quasi-coherent
sheaves respectively.

As we said we will apply \cite[Theorem 2.2]{hovey2} to get our
locally projective model structure. The adapted version of the
conditions of that theorem to our setting is the following:
\begin{enumerate}
\item The pairs $(\barP,\dgUo)$ and $(\dgP,\barUo)$ are cotorsion pairs,
\item Exact dg-locally projective complexes in $\Ch(\Qco({\bf P^1}(k)))$ are locally projective,
that is, $\dgP\cap {\mathcal E}=\barP$ where $\mathcal E$ is the
class of all exact complexes of quasi-coherent sheaves on ${\bf
P^1}(k)$.
\item The pairs $(\barP,\dgUo)$ and $(\dgP,\barUo)$ are complete.
\end{enumerate}

So let us prove each one of the previous conditions. Before that, we
need to make more accurate the statements made in \cite[Lemma
3.8(7),(8)]{Gill}.
\begin{lem}\label{mist}
Let $\mathcal A$ be an abelian category and $\Ch{\mathcal A})$ its
corresponding category of un\-boun\-ded complexes. Let $M,N$ be any
complexes of $\Ch(\mathcal A)$ and $C$ any object of $\mathcal A$.
Then there exist monomorphisms of abelian groups
$$0\to\Ext^1_{\mathcal A}(C,Z_nN) \to \Ext^1_{\Ch(\mathcal A)}(\underline{C}[-n],
N)$$ and $$0\to\Ext^1_{\mathcal A}(M_n/B_nM,C) \to
\Ext^1_{\Ch(\mathcal A)}(M,\underline{C}[-n])$$
\end{lem}
\begin{proof}
Let us see the first claim. Let $0\to Z_nN\to T\to C\to 0$ be any
extension of $\Ext^1_{\mathcal A}(C,Z_n N)$. Then we form the
pushout of the inclusion $Z_n N\to N_n$ and $Z_n\to T$ \DIAGV{60} {}
\n{} \n{0}\n{}\n{0}\nn {}\n{}\n{\sar}\n{}\n{\sar}\nn{0}
\n{\ear}\n{Z_n N}\n{\ear}\n{T}\n{\ear}\n{C}\n{\ear}\n{0}\nn
{}\n{}\n{\sar}\n{}\n{\sar}\n{}\n{\seql}\nn
{0}\n{\ear}\n{N_n}\n{\Ear{g}}\n{Q}\n{\ear}\n{C}\n{\ear}\n{0}\nn
{}\n{}\n{\sar}\n{}\n{\saR{h}}\nn
{}\n{}\n{B_{n+1}N}\n{\eeql}\n{B_{n+1}N}\nn
{}\n{}\n{\sar}\n{}\n{\sar}\nn {}\n{}\n{0}\n{}\n{0}\diag Then we have
the commutative \DIAGV{70} {0}
\n{\ear}\n{N_{n-1}}\n{\eeql}\n{N_{n-1}}\n{\ear}\n{0}\n{\ear}\n{0}\nn
{}\n{}\n{\Sar{\delta^{n-1}}}\n{}\n{\saR{g\circ
\delta^{n-1}}}\n{}\n{\sar}\n{}\nn
{0}\n{\ear}\n{N_n}\n{\Ear{g}}\n{Q}\n{\Ear{t}}\n{C}\n{\ear}\n{0}\nn
{}\n{}\n{\Sar{\delta^n}}\n{}\n{\Sar{h}}\n{}\n{\sar}\nn
{0}\n{\ear}\n{N_{n+1}}\n{\eeql}\n{N_{n+1}}\n{\ear}\n{0}\n{\ear}\n{0}\diag
Hence we have an extension $\xi=0\to N\to H\to \underline{C}[-n]\to
0$ in $\Ext^1_{\Ch(\mathcal A)}(\underline{C}[-n],N)$. This defines
a map from $$\Ext^1_{\mathcal A}(C,Z_n N)\to\Ext^1_{\Ch(\mathcal
A)}(\underline{C}[-n],N).$$It is clear that the map is going to be a
morphism of abelian groups with respect to the Baer sum of
extensions. Now if $\xi$ splits and $r:C\to Q$ is the corresponding
excision in the $n$th component of $\xi$, we follow by the
commutativity of the diagram that $h\circ r=0$, so ${\rm
Im}(r)\subseteq T$. Hence $0\to Z_n N\to T\to C\to 0$ splits. The
proof of the second monomorphism is dual.
\end{proof}

\begin{rem}
Monomorphisms of Lemma \ref{mist} are not isomorphisms in general.
For if we consider the category $\Ch(R)$ of unbounded complexes of
$R$-modules (here $R$ is any ring with identity) then if $P$ is any
projective $R$-module $\Ext^1_R(P,-)=0$, but
$\Ext^1_{\Ch(R)}(\underline{P}[-n], -)\neq 0$ because
$\underline{P}[-n]$ is a dg-projective complex, but is not
projective. The same holds for the second monomorphism by taking $C$
any injective $R$-module. We also point out that if $\mathcal{A}$
has enough injectives then we have the cotorsion pair
$$({\rm exact\ complexes},{\rm dg-injectives}) $$ and hence it
is easy to get the isomorphism
$$\Ext^1_{\mathcal A}(Z_{n+1} M,C)\cong \Ext^1_{\Ch(\mathcal
A)}(M,\underline{C}[-n]),$$ for every exact complex $M$. Dually if
$\mathcal{A}$ has enough projectives then we get the cotorsion pair
$$({\rm dg-projectives},{\rm exact\ complexes})$$and therefore if
$N$ is exact we have the isomorphism $\Ext^1_{\mathcal A}(C,Z_n
N)\cong \Ext^1_{\Ch(\mathcal A)}(\underline{C}[-n], N)$.
\end{rem}

\begin{prop}\label{soncot}
The pairs $(\barP,\dgUo)$ and $(\dgP,\barUo,)$ are cotorsion pairs,
\end{prop}
\begin{proof} Since $(\P,\Uo)$ is a cotorsion pair and $\P$ contains
a family of locally projective generators for $\Qco({\bf P^1}(k))$
the result follows from \cite[Proposition 3.6]{Gill} (with the
remark that the part of the proof of \cite[Proposition 3.6]{Gill}
involving \cite[Lemma 3.8(7),(8)]{Gill} can be replaced by Lemma
\ref{mist}).\end{proof}
\begin{prop}\label{dgp}
$\dgP\cap {\mathcal E}=\barP$ where $\mathcal E$ is the class of all
exact complexes of quasi-coherent sheaves on ${\bf P^1}(k)$.
\end{prop}
\begin{proof} By Lemma 3.10 of \cite{Gill} it only remains to prove
that $\dgP\cap {\mathcal E}\subseteq \barP$. By the results of
\cite[Section II.5]{Hartshorne} there exists right adjoint of the
restriction functor  $i^*[x]:\Qco({\bf P^1}(k))\to k[x]\Mod$ given
by $i^*(M_1\to M\fiz M_2)=M_1$. This is defined as
$i_*[x](N)=(N\hookrightarrow S^{-1}N\stackrel{id}{\fiz}S^{-1}N)$,
for every $k[x]$-module $N$. And there are an analogous pair of
adjoint functors $$(i^*[x^{-1}],i_*[x^{-1}])$$ and
$(i^*[x,x^{-1}],i_*[x,x^{-1}])$.

Now let $Y=Y_1\to Y_0\fiz Y_2$ be a complex in $\dgP\cap {\mathcal
E}$ (so $Y_1,Y_0$ and $Y_2$ are complexes of $k[x]$, $k[x^{-1}]$ and
$k[x,x^{-1}]$ modules respectively). To see that $Y$ is in $\barP$
we have to check that $Z_nY$ is a locally projective quasi-coherent
sheaf, for all $n\in \Z$, that is $Z_nY_1$, $Z_nY_0$ and $Z_nY_2$
are projective $k[x]$, $k[x^{-1}]$ and $k[x,x^{-1}]$ modules
respectively. By \cite[Proposition 2.3.7]{jr} if a complex of
modules (over $k[x]$, $k[x^{-1}]$ or $k[x,x^{-1}]$) is exact and
dg-projective then it is projective (so, in particular, $Z_nY_1$,
$Z_nY_0$ and $Z_nY_2$ will be projective modules). So we will be
done if we show that $Y_1$, $Y_0$ and $Y_2$ are exact and
dg-projective complexes of $k[x]$, $k[x^{-1}]$ and $k[x,x^{-1}]$
modules, respectively. We will do so for $Y_1$. The other cases are
similar. So let us assume that
$$Y_1=\cdots \to M^{-1}\to M^0\to M^1\to \cdots,$$ with $M^i\in
k[x]\Mod$, $i\in \Z$. Since $Y$ is an exact complex of
quasi-coherent sheaves $Y_1=i^*[x](Y)$ will be also exact. We see
that is dg-projective. So let $E$ be an exact complex of
$k[x]$-modules. We have to check that $\Hom_{\Ch(k[x])}(Y_1,E)$ is
exact. But, by the previous comments, there is an isomorphism
$$\Hom_{\Ch(k[x])}(Y_1,E)=\Hom_{\Ch(k[x])}(i^*[x](Y),E)\cong
\Hom_{\Ch(\Qco({\bf P^1}(k)))}(Y,i_*[x](E))$$ and since the functor
$i_*[x]$ preserves exactness, $i_*[x](E)$ will be an exact complex
of quasi-coherent sheaves on ${\bf P^1}(k)$. Since $Y\in \dgP$ if we
show that $i_*[x](E)\in \barUo$ we will be done. To see this we need
to check that $Z_ni_*[x](E)\in \Uo$, $\forall n\in \Z$. But
$Z_ni_*[x](E)=i_*[x](Z_nE)$. Hence $$\Ext^1_{\Qco({\bf
P^1}(k))}(\O(m),i_*[x](Z_nE))\cong
\Ext^1_{k[x]}(i^*[x](\O(m)),Z_nE)=0$$ (where the last equality
follows because $i^*[x](\O(m))=k[x]$).\end{proof}

To get that $(\dgP,\barUo)$ is a complete cotorsion pair we need the
following lemma.
\begin{lem}\label{reduc}
Let $Y$ be a complex in $\dgP$. Then $Y$ is a direct summand of a
direct transfinite extension of $\underline{\O(n)}[m]$'s.
\end{lem}
\begin{proof} We will prove that the cotorsion pair $(\dgP, \barUo)$
is cogenerated by the set $\mathcal{I}=\{\underline{\O(k)}[m]:
k,m\in \Z\}$. Then the result will follow reasoning as we did in the
proof of Theorem \ref{und}. It is easy to check that
$\mathcal{I}\subseteq \dgP$ for if $\underline{\O(k)}[m]_l\in \P,
\forall l\in \Z$ and for every exact complex $M\in \barUo$,
$Hom(\underline{\O(k)}[m],M)$ is the complex
$$\cdots\to \Hom(\O(k),M^l)\to
\Hom(\O(k),M^{l+1})\to\cdots$$which is obviously exact because
$Z_nM, B_nM\in \Uo$. So therefore $\mathcal{I}^{\perp}\supseteq
(\dgP)^{\perp}=\barUo$. We now prove the converse: let $N\in
\mathcal{I}^{\perp}$. We have to see that $N$ is exact and that
$Z_nN\in \Uo$. We prove that $N$ is exact. It is clear that this is
equivalent to that each morphism $\underline{\O(m)}[n]\to N$ can be
extended to $\overline{\O(m)}[n]\to N$, for every $m,n\in \Z$. But
this follows from the short exact sequence
$$0\to \underline{\O(m)}[n]\to \overline{\O(m)}[n] \to \underline{\O(m)}[n+1] \to 0 $$
and since $\Ext^1(\underline{\O(m)}[n+1],N)=0$. Let us see the last
claim. Since $\U=\{\O(m):\ m\in \Z\}$ cogenerates the cotorsion pair
$(\P,\Uo)$ we only need to prove that $\Ext^1_{\Qco({\bf
P^1}(k))}(\O(m),Z_n N)=0$, $\forall m,n\in \Z$. By Lemma \ref{mist}
we have a monomorphism of abelian groups $$0\to \Ext^1_{\Qco({\bf
P^1}(k))}(\O(m),Z_n N)\to \Ext^1_{\Ch(\Qco({\bf
P^1}(k)))}(\underline{\O(m)}[-n],N)$$ and since the last is equal to
0 we get that $Z_n N\in \Uo$.
\end{proof}

\begin{prop}\label{sufinye}
The cotorsion pair $(\dgP,\barUo)$ of complexes of quasi-coherent
sheaves on ${\bf P^1}(k)$ is complete.
\end{prop}
\begin{proof} By Lemma \ref{reduc} the pair $(\dgP,\barUo)$ is a cotorsion
pair cogenerated by a set and $\dgP$ contains a family of generators
of $\Ch(\Qco({\bf P^1}(k))$. Then by \cite[Theorem 2.6]{Estr} we get
that $(\dgP,\barUo)$ is a complete cotorsion pair.\end{proof}
\begin{cor}
Let $\mathcal{E}$ be the class of exact complexes of quasi-coherent
sheaves on ${\bf P^1}(k)$, then $\barUo=\dgUo\cap \mathcal{E}$.
\end{cor}
\begin{proof} By propositions \ref{soncot} and \ref{sufinye} the pair
$(\dgP,\barUo)$ is a cotorsion pair with enough injectives. By
Proposition \ref{dgp}, $\dgP\cap \mathcal{E}=\barP$, so by
\cite[Lemma 3.14 (a)]{Gill} we get that we claim. \end{proof}

We finish this section by proving that $(\barP,\dgUo)$ is also
complete. We need the following Lemma.
\begin{lem}\label{eskap}
The class $\P$ of all locally projective quasi-coherent sheaves on
${\bf P^1}(k)$ is a Kaplansky class.
\end{lem}
\begin{proof}
Let $P\in \P$ be a locally projective quasi-coherent sheaf. By
Theorem \ref{p1} we can write $P=\li_{\alpha<\lambda}S_{\alpha}$,
with $\{S_{\alpha}:\ \alpha<\lambda\}$ a direct transfinite system
of countable generated quasi-coherent sheaves on ${\bf P^1}(k)$. Let
$\aleph\geq \omega,|k|$ be a regular cardinal and let $0\neq
X\subseteq P$ where $|X|\leq \aleph$. For every element $x\in X$ let
us pick $j_x<\lambda$ such that $x\in S_{j_x}$. Let $\gamma$ be the
supremmum of such $j_x$, $x\in X$ and let us take
$S=\li_{\beta<\gamma}S_{\beta}$. It is clear that $|S|\leq \aleph$
and that $S\in \P$. Let us see that $P/S\in \P$. Since direct limits
in $\Qco({\bf P^1}(k))$ are computed componentwise, if we call
$S=S_1\to S_0\fiz S_2$ and $P=P_1\to P_0\fiz P_2$, we get that
$S_1=\li_{\beta\leq \gamma}S_{\beta}^1$ (where
$S_{\beta}=S_{\beta}^1\to S_{\beta}^0\fiz S_{\beta}^2$) and $P_1$
are direct transfinite extensions of countably generated projective
$k[x]$-modules, so hence direct sums of countably generated
projective $k[x]$-modules and $S_1$ is a direct summand of $P_1$.
Therefore $P_1/S_1$ will be a projective $k[x]$-module. The same
reasoning applies to $P_0/S_0$ and $P_2/S_2$ to get that $P/S$ is a
locally projective quasi-coherent sheaf.
\end{proof}

\begin{thm}
The cotorsion pair $(\barP,\dgUo)$ is complete.
\end{thm}

\begin{proof} We will make the proof in several steps. First we will
use Lemma \ref{eskap} to see that the pair $(\barP,\dgUo)$ is
cogenerated by a set. Then we appeal to \cite[Theorem 2.6]{Estr} to
get that the cotorsion pair is complete (note that the set
$\{\overline{\O(k)}[m]: k,m\in \Z\}\subseteq \barP$ also cogenerates
$\Ch(\Qco({\bf P^1}(k)))$). To see that the pair $(\barP,\dgUo)$ is
cogenerated by a set we need to show the following: let $P$ be any
exact complex in $\barP$, $x\in P$ and let us fix a regular cardinal
$\aleph\geq \omega,|k|$. We will prove that there exists an exact
subcomplex $L$ of $F$ such that $L,F/L\in \barP$ and $|L|\leq
\aleph$. Since the class $\barP$ is closed under extensions and
direct limits the previous says that we can write every complex in
$\barP$ as the direct union of a continuous chain of subcomplexes in
$\barP$ with cardinality less than or equal to $\aleph$. Then if $T$
is a set of representatives of complexes $L$ in $\barP$ with
$|L|\leq \aleph$, we get by \cite[Lemma 1]{EkTrl} that the pair
$(\barP,\dgUo)$ is cogenerated by a set.

So let us start with the proof. We fix some notation: let us denote
by $G=\oplus_{n\in\Z}\O(n)$ the generator of $\Qco({\bf P^1}(k))$.
For a given $x\in F^i$ there exists a certain $m\in \Z$ and a map
$\O(m)\to F^i$ sending 1 to $x$. we use the notation $Gx$ to denote
the image of map $\O(m)\hookrightarrow G\to F^i$. Let us suppose
(without loss of generality) that $k=0$ and $x\in F^0$. Consider
then the exact complex
$$(S1)\hskip 1cm \cdots \rightarrow
A^{-2}_1\stackrel{\delta^{-2}}{\rightarrow} A^{-1}_1
\stackrel{\delta^{-1}}{\rightarrow} Gx
\stackrel{\delta^{0}}{\rightarrow} \delta^0(Gx)
\stackrel{\delta^1}{\rightarrow} 0$$ where $A^{-i}_1$ is a
quasi-coherent subsheaf of $F^{-i}$ constructed as follows:
$|Gx|\leq \aleph$ since $|G|\leq\aleph$, so we can find
$A^{-1}_1\leq F^{-1}$ such that $|A^{-1}_1|\leq\aleph$ and
$\delta^{-1} (A^{-1}_1)=\ker (\delta^0|_{Rx})$. Then $A^{-2}_1\leq
F^{-2}$, $|A^{-2}_1|\leq\aleph$, and $\delta^{-2}(A^{-2}_1)=\ker
(\delta^{-1} |_{A^{-1}_1})$, and we repeat the argument.

Now $\ker(\delta^0|_{Gx})\leq\ker \delta^0$, so we know by Lemma
\ref{eskap} that $\ker (\delta^0|_{Gx})$ can be embedded into a
locally projective quasi-coherent subsheaf $S^0_2$ of $\ker
\delta^0$. Since $|\ker (\delta^0|_{Gx})|\leq\aleph$ we see by Lemma
\ref{eskap} that $S^0_2$ can be chosen in such a way that
$|S^0_2|\leq\aleph$. Then consider the exact complex
$$(S2)\hskip 1cm \cdots \rightarrow
A^{-2}_2\stackrel{\delta^{-2}}{\rightarrow} A^{-1}_2
\stackrel{\delta^{-1}}{\rightarrow} Gx +S^0_2
\stackrel{\delta^{0}}{\rightarrow} \delta^0(Gx)
\stackrel{\delta^1}{\rightarrow} 0$$ where $A^{-i}_2$ are taken as
above. It is clear that $\ker (\delta^0|_{Gx+S^0_2}) =S^0_2$, which
is a locally projective quasi-coherent subsheaf of $\ker \delta^0$,
and that $|Gx+S^0_2|\leq\aleph +\aleph =\aleph$.

Observe now that $\delta^0(Gx)\subseteq \ker \delta^1$, so we can
embed $\delta^0(Gx)$ into a locally projective quasi-coherent
subsheaf $S^1_3$ of $\ker\delta^1$ in such a way that $|S^1_3|\leq
\aleph$ ($|\delta^0(Gx)|\leq\aleph$), and then take the exact
complex
$$(S3)\hskip 1cm \cdots \rightarrow
A^{-2}_3\stackrel{\delta^{-2}}{\rightarrow} A^{-1}_3
\stackrel{\delta^{-1}}{\rightarrow} A^0_3
\stackrel{\delta^{0}}{\rightarrow} S^1_3
\stackrel{\delta^1}{\rightarrow} 0.$$ We see again that $\ker
(\delta|_{S^1_3}) =S^1_3$, which is a quasi-coherent subsheaf of
$\ker\delta^1$ in $\P$.

We turn over and find $S^0_4\leq\ker\delta^0$ locally projective
with $|S^0_4|\leq\aleph$ and $S^0_4\supseteq\ker
(\delta^0|_{A^0_3})$, and then construct $A^{-i}_4\leq F^{-i}$
($|A^{-i}_4|\leq\aleph\ \forall i$) such that $$(S4)\hskip 1cm
\cdots \rightarrow A^{-2}_4\stackrel{\delta^{-2}}{\rightarrow}
A^{-1}_4 \stackrel{\delta^{-1}}{\rightarrow} A^0_3 +S^0_4
\stackrel{\delta^{0}}{\rightarrow} S^1_3
\stackrel{\delta^1}{\rightarrow} 0$$ is exact. Once more $\ker
(\delta^0|_{A^0_3+S^0_4})= S^0_4\leq\ker\delta^0$ is a locally
projective quasi-coherent subsheaf. Then find
$S^{-1}_5\leq\ker\delta^{-1}$ locally projective with
$|S^{-1}_5|\leq\aleph$, $\ker (\delta^{-1}|_{A^{-1}_4})\subseteq
S^{-1}_5$, and consider the exact complex $$(S5)\hskip 1cm \cdots
\rightarrow A^{-2}_5\stackrel{\delta^{-2}}{\rightarrow} A^{-1}_4+
S^{-1}_5 \stackrel{\delta^{-1}}{\rightarrow} A^0_3 +S^0_4
\stackrel{\delta^{0}}{\rightarrow} S^1_3
\stackrel{\delta^1}{\rightarrow} 0,$$ in which
$\ker(\delta^{-1}|_{A^{-1}_4+S^{-1}_5})
=S^{-1}_5\leq\ker\delta^{-1}$ pure.

The next step is to find $S^{-2}_6\leq\ker\delta^{-2}$ locally
projective such that $|S^{-2}_6|\leq\aleph$ and that
$\ker(\delta^{-2}|_{A^{-2}_5})\subseteq S^{-2}_6$, and then consider
the exact complex $$(S6)\hskip 1cm \cdots \rightarrow
A^{-3}_6\stackrel{\delta^{-3}}{\rightarrow}
A^{-2}_5+S^{-2}_6\stackrel{\delta^{-2}}{\rightarrow} A^{-1}_4+
S^{-1}_5 \stackrel{\delta^{-1}}{\rightarrow} A^0_3 +S^0_4
\stackrel{\delta^{0}}{\rightarrow} S^1_3
\stackrel{\delta^1}{\rightarrow} 0$$ in which
$\ker(\delta^{-2}|_{A^{-2}_5+S^{-2}_6})=S^{-2}_6 \subseteq
\ker\delta^{-2}$ locally projective.

Therefore we prove by induction that for any $n\geq 4$ we can
construct an exact complex $$(Sn)\cdots
\stackrel{\delta^{-n+2}}{\rightarrow} A^{-n+3}_n
\stackrel{\delta^{-n+3}}{\rightarrow} T^{-n+4}_n
\stackrel{\delta^{-n+4}}{\rightarrow} T^{-n+5}_n \rightarrow \cdots
\stackrel{\delta^{-1}}{\rightarrow}
T^0_n\stackrel{\delta^0}{\rightarrow} T^1_n
\stackrel{\delta^1}{\rightarrow} 0$$ such that $\ker
(\delta^{-n+j}|_{T^{-n+j}_n})$ is a locally projective
quasi-coherent subsheaf of $\ker\delta^{-n+j}$ $\forall j\geq 4$ and
that all the terms have cardinality less than or equal to $\aleph$.

If we take the direct limit $L={\displaystyle \lim_{\rightarrow}}
(Sn)$ with $n\in\Natur$, we see that the complex $L$ is exact and
$\ker(\delta^i|_{L^i})$ is a locally projective quasi-coherent
subsheaf of $\ker\delta^i$ $\forall i\leq 1$. Furthermore
$|L^i|\leq\aleph_0\cdot\aleph= \aleph$ for any $i\leq 1$, so
$|L|\leq\aleph$. We finally consider the complex $L$ to be
$$L=\cdots \rightarrow L^i\stackrel{\delta^i}{\rightarrow} L^{i+1}
\stackrel{\delta^{i+1}}{\rightarrow} \cdots
\stackrel{\delta^{-1}}{\rightarrow}
L^0\stackrel{\delta^0}{\rightarrow} L^1
\stackrel{\delta^1}{\rightarrow} 0 \stackrel{\delta^2}{\rightarrow}
0\cdots,$$ which is a subcomplex of $F$, $x\in L^0$, and
$\ker(\delta^i|_{L^i})$ is a locally projective quasi-coherent
subsheaf of $\ker\delta^i$ $\forall i\in\Z$ and so
$\ker(\delta^i|_{L^i})$. Therefore the complex $L$ is a subcomplex
in $\barP$ of $F$ and of course $|L|\leq\aleph$.

To finish the proof we only have to argue that
$F/L=(F^i/L^i,\overline{\delta}^i)$ is in $\barP$. An easy
computation shows that
$\ker\overline{\delta}^i=\ker(\delta^i)/\ker(\delta^i|_{L^i})$, but
by construction $\ker(\delta^i|_{L^i})$ is a locally projective
quasi-coherent subsheaf of $\ker\delta^i$ $\forall i\in\Z$, so
$\ker(\delta^i)/\ker(\delta^i|_{L^i})$ is locally projective for all
$i\in \Z$. Of course $F/L$ is exact since both $F$ and $L$ are
exact, so $F/L$ is in $\barP$.
\end{proof}

\section{The monoidal locally projective model structure on $\Ch(\Qco({\bf
P^1}(k)))$}\label{fin}

With the results of the previous section, we are in position to
impose a locally projective model structure on $\Ch(\Qco({\bf
P^1}(k)))$.

\begin{thm}
There is a model structure in $\Ch(\Qco({\bf P^1}(k))$ such that
$\dgP$ is the class of cofibrant objects, $\dgUo$ is the class of
fibrant objects and the exact complexes are the trivial objects.
\end{thm}
\begin{proof} This follows from \cite[Theorem 2.2]{hovey2} taking
$\mathcal{C}=\dgP$, $\mathcal F=\dgUo$ and
$\mathcal{W}=\mathcal{E}$, the class of all exact complexes of
quasi-coherent sheaves. \end{proof}

\bigskip
Now we will prove that the previous model structure is compatible
with the graded tensor product on $\Ch(\Qco({\bf P^1}(k))$. We
recall that for a given two complexes of quasi-coherent sheaves $M$
and $N$, the tensor product $M\otimes N$ is a complex of abelian
groups with $(M\otimes N)_m=\oplus_{t\in \Z}M_t\otimes_{\Qco({\bf
P^1}(k))}N_{m-t}$ and $$\delta=\delta_M^t\otimes
id_N+(-1)^tid_M\otimes \delta^{m-t}_N,$$ for all $m,t\in \Z$.

The previous tensor product becomes $\Ch(\Qco({\bf P^1}(k))$ into a
monoidal category. To see that the structure is closed we appeal to
the natural embedding $\Qco({\bf P^1})(k)\hookrightarrow \O_{{\bf
P^1}(k)}\Mod$, where $\O_{{\bf P^1}(k)}\Mod$ is the category of
sheaves of $\O_{{\bf P^1}(k)}$-modules. Since this embedding
preserves direct limits, it will have a right adjoint functor $Q:
\O_{{\bf P^1}(k)}\Mod\to \Qco({\bf P^1})(k)$. This functor $Q$
allows to show that $\Qco({\bf P^1})(k)$ is a closed symmetric
monoidal category (the closed structure is given by applying $Q$
after the internal Hom functor of $\O_{{\bf P^1}(k)}\Mod$). This
structure extends to $\Ch(\Qco({\bf P^1}(k)))$ becoming it into a
closed symmetric monoidal category

As it is pointed in \cite[pg. 9]{Hov2} it would be desirable to get
a model structure in the category $\Ch(\Qco({\bf P^1})(k))$
compatible with the closed symmetric monoidal structure (in the
sense of \cite[Chapter 4]{Hov}. Our locally projective model
structure certainly is (Theorem \ref{esmonoi}). We remark that for
the case where $X$ is a quasi-compact and quasi-separated scheme the
category $\Qco(X)$ has enough flat objects, so by using the results
of \cite{Est}, a modified argument to that of \cite{Gill} allows one
to impose a flat model structure in $\Ch(\Qco(X))$ which will be
compatible with the tensor product. However, since quasi-coherent
sheaves play the role of the modules in categories of sheaves and is
known that there exists a monoidal projective model structure in
$\Ch(R)$ (whenever $R$ is any commutative ring) it seems natural to
conjecture that there is analogous locally projective monoidal model
structure for quasi-coherent sheaves, at least for sufficiently nice
schemes (closed subschemes of ${\bf P}^m(k)$). So our result is a
first step in addressing this problem.

\begin{thm}\label{esmonoi}
The induced model structure on $\Ch(\Qco({\bf P^1}(k)))$ by the
cotorsion pair $(\P,\Uo) $ is compatible with the graded tensor
product given above.
\end{thm}
\begin{proof} Let us check that the conditions of \cite[Theorem
7.2]{hovey2} holds in this situation. Notice that using the notation
of that Theorem in our situation ${\mathcal P}$ is the class of all
short exact sequences and ${\mathcal W}$ the class of exact
complexes. So we will check that
\begin{itemize}
\item[$i)$] Every monomorphism of complexes with cokernel
a dg-locally projective complex is a pure injection in each degree.

\item[$ii)$] If $X$ and $Y$ are dg-locally projective
complexes then $X\otimes Y$ is also dg-locally projective.
\item[$iii)$] If $X,Y$ are dg-locally projective complexes and
$Y$ is exact then $X\otimes Y\in \barP$.
\item[$iv)$] The complex with the direct sum of $\O(m)$'s in
one component and 0 in the rest is dg-locally projective.

\end{itemize}
Conditions $i)$ and $iv)$ follows immediately from the definitions
(since a dg-locally projective complex is a flat quasi-coherent
sheaf componentwise). So let us see condition $ii)$. By Lemma
\ref{reduc} it suffices to prove the statement for
$\underline{\O(m)}[n]$, $n,m\in\Z$. But in this case we have
$$\underline{\O(m)}[n_1]\otimes \underline{\O(m')}[n_2]\cong
\underline{\O(m+m')}[n_1+n_2]$$ so is again of this form. Finally
let us check condition $iii)$. By $ii)$, $X\otimes Y$ is in $\dgP$
and since $Y$ is exact $X\otimes Y$ will be also exact. But then by
Proposition \ref{dgp} we get that $X\otimes Y\in \barP$.
\end{proof}

\section{Derived extension functors by using locally projective
resolutions}

We finish the paper by applying the previous monoidal locally
projective model structure to get an alternate way of computing
right derived functors of quasi-coherent sheaves on ${\bf P^1}(k)$
by using locally projective resolutions of quasi-coherent sheaves.

Derived extension functors can be defined from a Quillen's model
structure on $$\Ch(\Qco({\bf P^1}(k)))$$ from the equation
$$\Ext^n_{\Qco({\bf P^1}(k))}(M,N)=\Hom_{\Qco({\bf
P^1}(k))}(Q_M,R_N)/\sim$$where $Q_M$ is a cofibrant replacement of
$\underline{M}$ and $R_N$ is a cofibrant replacement of
$\underline{N}[n]$. Now in \cite[Lemma 5.3]{Gill} is shown that a
cofibrant replacement of $\underline{M}$ is the deleted complex $
{\bf P}_{\bullet}$ of an exact complex $\cdots \to P_2\to P_1\to
P_0\to M\to 0$ where $P_n$ are locally projective quasi-coherent
sheaves and each cycle quasi-coherent sheaf is in $\Uo$. Dually the
deleted complex ${\bf U}^{\bullet} $ of an exact complex $0\to N\to
U_0\to U_1\to \cdots$ where $U_n$ is in $\Uo$ and each cycle is
locally projective, is a fibrant replacement of $\underline{N}$.

Now noting that $$\Hom_{\Qco({\bf P^1}(k))}({\bf P}_{\bullet},{\bf
U}^{\bullet}[n])/\sim=H^n(Hom({\bf P}_{\bullet},{\bf
U}^{\bullet}))$$ we follow that derived extension functors of
quasi-coherent sheaves can be computed as the homology of the
$Hom$-complex:
$$\Ext^n_{\Qco({\bf
P^1}(k))}(M,N)=H^n(Hom({\bf P}_{\bullet},{\bf U}^{\bullet}))$$ where
${\bf P}_{\bullet}=\cdots \to P_2\to P_1\to P_0\to 0$ and ${\bf
U}^{\bullet}= 0\to U_0\to U_1\to U_2\to \cdots$ are the
corresponding deleted complexes of $M$ and $N$.

\end{document}